\documentclass[a4paper,12pt]{article}

\usepackage[utf8]{inputenc} 
\usepackage{amsthm,amsmath,amscd,mathtools}
\usepackage[vmargin=2cm,hmargin=2cm,headheight=14.5pt,top=2cm,headsep=.5cm]{geometry}
\usepackage{hyperref}
\usepackage[autobold]{mathfixs}
\usepackage{stackrel}
\usepackage[charter]{mathdesign}
\usepackage{bm}
\usepackage[singlelinecheck=false]{caption}
\usepackage{xcolor}
\usepackage{graphicx}
\DeclareMathSizes{12}{12}{8}{6}

\usepackage{cite}
\usepackage{enumitem,subfigure}
\usepackage{pgfplots}
\setlist[itemize,2]{label=$\centerdot$}
\setlist[itemize,3]{label=$\triangle$}
\usepackage{tikz-cd}
\usepackage{tikz}
\usepackage{graphicx} 

\newtheoremstyle{ptheorem}{1em}{0em}{\itshape}{}{\bfseries}{.}{.5em}{\thmname{#1}\thmnumber{
		#2}\thmnote{ (\hspace{-1sp}{#3})}}

\theoremstyle{ptheorem}

\newtheorem{thm}{Theorem}[section]
\newtheorem{pro}[thm]{Proposition}
\newtheorem{lem}[thm]{Lemma}
\newtheorem{cor}[thm]{Corollary}

\newtheoremstyle{hdef}{1em}{0em}{}{}{\bfseries}{.}{.5em}{\thmname{#1}\thmnumber{
		#2}\thmnote{ (\hspace{-.01pt}{#3})}}
\theoremstyle{hdef}

\newtheorem{dfn}[thm]{Definition}
\newtheorem{rem}[thm]{Remark}

\numberwithin{equation}{section}
\numberwithin{figure}{section}

\newcommand{\bN}{{\mathbb N}}

\renewcommand{\a}{\alpha}
\renewcommand{\b}{\beta}

\renewcommand{\phi}{\varphi}

\renewcommand{\(}{\left(}
\renewcommand{\)}{\right)}

\newcommand{\bs}{\backslash}

\newcommand{\olb}[1]{%
	\vbox{\offinterlineskip\ialign{\hfil##\hfil\cr
			$\rotatebox[origin=c]{90}{$]$}$\cr\noalign{\kern-.45ex}{$#1$}\cr}}}

\newcommand{\noop}[1]{}

\usepackage{stmaryrd}

\newcommand\restr[2]{{
 \left.\kern-\nulldelimiterspace 
 #1 
 \littletaller 
 \right|_{#2} 
 }}
\newcommand{\littletaller}{\mathchoice{\vphantom{\big|}}{}{}{}}

\allowdisplaybreaks

\begin{document}

\title{On the approximation properties of Stieltjes polynomials}
\author{V\'ictor Cora and F. Adri\'an F. Tojo}

\date{}

\author{Víctor Cora$^{*}$\\
	\normalsize e-mail: victor.cora.calvo@usc.es\\
	F. Adri\'an F. Tojo$^{*\dagger}$ \\
	\normalsize e-mail: fernandoadrian.fernandez@usc.es\\ \normalsize \emph{$^{*}$Departamento de Estatística, Análise Matemática e Optimización},\\ \normalsize \emph{Universidade de Santiago de Compostela, 15782, Facultade de Matemáticas, Santiago, Spain.}\\\normalsize
	\\\normalsize \emph{$^{\dagger}$CITMAga, 15782, Santiago de Compostela, Spain}}

\date{}

\maketitle

\medbreak

\begin{abstract}
We introduce and study the approximation properties of $g$--polynomials, defined as linear combinations of iterated Stieltjes integrals of a constant function. Focusing on the case where the derivator $g$ has finitely many discontinuities, we prove that the space of $g$--polynomials is dense in the space of uniformly $g$--continuous functions. This result establishes a Weierstrass-type approximation theorem within the framework of Stieltjes calculus. The characterization of the closure of the space of $g$--polynomials in the general case, where the derivator may exhibit more complex behavior, remains an open and challenging problem, which we briefly discuss.
\end{abstract}


\medbreak


\section{Introduction}

Stieltjes differential equations have been widely studied since they naturally model phenomena with impulses or that have latent states
---see, for instance, \cite{LopezPousoMarquez2019}. The more theoretical and analytical tools we can provide for Stieltjes calculus, the easier these differential equations become to study. This work aims to help in that direction.

In classical analysis, the Weierstrass Approximation Theorem asserts that every continuous function on a compact interval $[a,b]$ may be uniformly approximated by polynomials \cite{1885}. This is a classical, fundamental result that underpins a vast array of methods from Fourier analysis to numerical schemes. We aim to obtain an analogous result for Stieltjes Calculus.

In this setting, the analogue of a polynomial is the \emph{$g$--polynomial}, a linear combination of iterated Stieltjes integrals of a constant function. This concept is widely studied in \cite{Cora2023}. The main contribution of this paper is to establish that \emph{whenever the derivator has only finitely many points of discontinuity}, the space of $g$--polynomials, ${\rm P}_g$, remains dense in $\mathrm{UC}_g([a,b])$, the space of uniformly $g$--continuous functions. Although we are not proving this result in full generality, since derivators may have infinitely many discontinuities, in practice, all real world models involve only finitely many jump points, so that solutions are computable.

When the derivator admits jump discontinuities the topology induced by the pseudometric disconnects the interval at the jump points and the algebraic structure of ${\rm P}_g$ is no longer closed under multiplication. These obstacles make the problem very challenging, since the usual strategies used to tackle polynomial approximation are not easily adaptable to the Stieltjes setting.

The paper is organized as follows. Section~\ref{seccionpre} provides the necessary background on Stieltjes calculus and Lebesgue-Stieltjes integration, including the definition and key properties of $g$--monomials and $g$--polynomials, with particular reference to~\cite{Cora2023}. In Section~\ref{seccionrutas}, we summarize the various approaches considered in our attempts to prove the general case and briefly comment on the challenges encountered. We also reduce the problem to the study of a suitable family of functions and show that this can be approached through the asymptotic behavior of Gram matrices. Since our method for proving the density of $\mathrm{P}_g$ cannot be generalized to the case of infinitely many discontinuities, we leave it as a potential tool for future work. Section~\ref{seccionaux} presents a collection of auxiliary results required for the main theorem, which constitute the core of the theoretical development. In Section~\ref{pruebaseccion}, we establish a Weierstrass-type approximation theorem within the framework of Stieltjes calculus for derivators with finitely many discontinuities. Finally, in Section~\ref{secfinal}, we extend this result to derivators with only a jump part (i.e., $g^C = 0$), under a suitable topological condition on the infinite set of discontinuities. Concluding remarks are given in Section~\ref{finalcooncl}.

\section{Preliminaries}
\label{seccionpre}

Let $[a,b]\subset \mathbb R$ be an interval. Let $\mathbb F\in\{\mathbb R,\mathbb C\}$. We will say that a nondecreasing and left-continuous function $g:[a,b]\to \mathbb R$  is a \emph{derivator}.
This derivator defines an outer measure that restricts to a measure on a Borel $\sigma$-algebra over the reals following Carath\'eodory's construction; for more details see \cite{Pouso2015}. Let us denote this measure by $\mu_g$.
The integral defined through this measure will be called the \emph{Lebesgue--Stieltjes integral}.

Given an interval $[c,d)\subset [a,b]$, we have $\mu_g([c,d))=g(d)-g(c)$. The pseudometric $d_g(x,y)=|g(x)-g(y)|$, $x,y\in[a,b]$, allows defining a topology by taking balls:
\[
B_g(x,\delta)=\{y\in[a,b]:|g(x)-g(y)|<\delta\}.
\]
The open sets are arbitrary unions of such balls. Let us denote this topology by $\tau_g$.
\begin{dfn}
 We say that a function $f:[a,b]\to\mathbb F$ is \emph{$g$--continuous}
 at $x\in[a,b]$ if and only if
 \[
 \forall \varepsilon >0 \,\,\exists \delta(x,\varepsilon) >0 \,\text{ such that if } y\in B_g(x,\delta)\Rightarrow |f(x)-f(y)|<\varepsilon.
 \]
 A function is $g$--continuous if it is $g$--continuous at every point of its domain. We say a function $f$ is \emph{uniformly $g$--continuous} if for every $\varepsilon>0$ there exists $\delta(\varepsilon)>0$ such that
 \[
 \forall x,y\in[a,b] \text{ such that } |g(x)-g(y)|<\delta(\varepsilon) \Rightarrow |f(x)-f(y)|<\varepsilon.
 \]
\end{dfn}
Define
\[
\Delta g(x)=g(x^+)-g(x),
\]
for $x\in[a,b)$. Here, $g(x^+)$ denotes the right-hand limit of $g$ at $x$. Consider the following sets:
\[
D_g=\{x\in[a,b):\Delta g(x)>0\}
\]
and
\[
C_g=\{x\in[a,b]:\exists \delta>0 \text{ such that } g \text{ is constant in } (x-\delta,x+\delta)\cap [a,b]\}.
\]
Note that $C_g$ is a usual open subset of $[a,b]$, $D_g$ is the set of discontinuity points of $g$ and, hence, at most countable.

The derivator $g$ has been used to define the Stieltjes derivative of a real-valued function \cite{Pouso2015}, hence its name.

\begin{dfn}[{\cite[Definition 3.7]{Fernandez2022}}]
	We define the \emph{Stieltjes derivative}\index{Stieltjes derivative}, or \emph{$g$--derivative}\index{$g$--derivative}, of a function $f:{\mathbb R}\to {\mathbb R}$ at a point $t\in {\mathbb R}$ as
	\[
	f'_g(t)=
	\begin{dcases}
		\displaystyle \lim_{s \to t}\frac{f(s)-f(t)}{g(s)-g(t)}, & t\not\in D_{g}\cup C_g,\\
		\displaystyle\lim_{s\to t^+}\frac{f(s)-f(t)}{g(s)-g(t)}, & t\in D_{g},\\
		\displaystyle\lim_{s\to b_n^+}\frac{f(s)-f(b_n)}{g(s)-g(b_n)}, & t\in C_{g},\ t\in(a_n,b_n),
	\end{dcases}
	\]
	where each $(a_n,b_n)$ is a connected component of $C_g$, provided the corresponding limits can be considered and exist. In that case, we say that $f$ is \emph{$g$--differentiable at $t$}.
\end{dfn}

The derivator $g$ can be decomposed into its continuous part and its jump part. Define, for $x\in[a,b]$,
\[
g^B(x)=\sum_{y\in[a,x)\cap D_g}\Delta g(y)
\]
and $g^C(x)=g(x)-g^B(x)$. Both $g^C$ and $g^B$ are left-continuous and nondecreasing, and therefore, derivators. In particular $g^C$ is continuous at every point.

Both notions of integral and Stieltjes derivative are consistent, and we have the following result.
\begin{thm}[{\cite[Theorem 5.4]{Pouso2015}}]
 \label{tf}
 Let $F:[a,b] \to \mathbb R$. The following conditions are equivalent:
 \begin{enumerate}
 \item[(1)] The function $F$ is \emph{$g$--absolutely continuous}: for every $\varepsilon>0$, there exists some $\delta>0$ such that, for any family $\{(a_n,b_n)\}_{n=1}^{m}$ of pairwise disjoint open subintervals of $[a,b]$, the inequality
 \[
 \sum_{n=1}^{m}(g(b_n)-g(a_n))<\delta
 \]
 implies
 \[
 \sum_{n=1}^m|F(b_n)-F(a_n)|<\varepsilon.
 \]
 \item[(2)] The function $F$ fulfills the following properties:
 \begin{enumerate}
 \item[(a)] there exists $F'_g(t)$ for $\mu_g$--almost all $t\in [a,b)$ (i.e., for all $t$ except on a set of $\mu_g$ measure zero);
 \item[(b)] $F'_g \in \operatorname{L}^1_g([a,b))$, the set of Lebesgue--Stieltjes integrable functions with respect to $\mu_g$;
 \item[(c)] for each $t \in [a,b]$, we have
 \begin{equation*}
 F(t)=F(a)+\int_{[a,t)} F'_g(s) \operatorname*{d}\mu_g.
 \end{equation*}
 \end{enumerate}
 \end{enumerate}
 \end{thm}

In \cite{Cora2023}, the authors define the notion of monomial in the context of Stieltjes differentiation through recurrent integrals of constant functions.
\begin{dfn}[{\cite[Definition 3.1]{Cora2023}}]
 \label{gmonomio}
Let $x_0\in[a,b]$. Define
$g_{x_0,0}(x)=1$ for all $x\in[a,b]$. Recursively, define the \emph{$g$--monomial of degree $n$ and center $x_0$} as
\[
g_{x_0,n}(x)=\begin{dcases}
 n\int_{[x_0,x)} g_{x_0,n-1}\operatorname*{d}\mu_g, & x\geq x_0,\\
 -n\int_{[x,x_0)} g_{x_0,n-1}\operatorname*{d}\mu_g, & x<x_0,
\end{dcases}
\]
for $x\in[a,b]$.
\end{dfn}
 In \cite[Section 3]{Cora2023}, they study the $g$--monomials extensively. The $g$--monomials are infinitely $g$--differentiable functions. In \cite[Section 4]{Cora2023}, these $g$--monomials are used to define power series and thereby $g$--analytic functions.

By Theorem~\ref{tf}, the $g$--monomials are $g$--absolutely continuous functions and, therefore, uniformly $g$--continuous.

Following \cite{Cora2023}, we define the  \emph{$g$--polynomials} as linear combinations of $g$--monomials and we will denote by ${\rm P}_g$ the set of $g$--polynomials. From \cite[Proposition 3.14]{Cora2023}, it follows that a $g$--monomial centered at a point $x_1\in[a,b]$ is a linear combination of $g$--monomials centered at
$x_0$. Hence, for a fixed arbitrary $x_0\in[a,b]$, the set $\{g_{x_0,n}\}_{n=0}^\infty$ is a generating set for ${\rm P}_g$. The following result is of particular interest.

\begin{thm}[{\cite[Theorem 3.22]{Cora2023}}]
 \label{formulasmonomios}
 Let $g:[a,b]\to\mathbb R$ be a derivator and fix $x_0\in[a,b]$. For $n\in\mathbb N$ and $x\in[a,b]$,
 \[
 g_{x_0,n}(x)=\sum_{k=0}^n {n\choose k} g^C_{x_0,k}(x)g^B_{x_0,n-k}(x).
 \]
\end{thm}

The above result illustrates the relation between $g$--monomials with $g^C$--monomials and $g^B$--monomials. The above will be a key factor in the proof of the Weierstrass Approximation Theorem for Stieltjes Calculus.

Additionally, gaining a better understanding of the $g^B$ and $g^C$--monomials provides valuable insights into $g$--monomials in general. We present a brief summary of some results proven in \cite{Cora2023}.

\begin{pro}[{\cite[Lemma 3.5, Propositions 3.15-3.16]{Cora2023}}]
 \label{propiedadesmonomios}
 Let $g:[a,b]\to \mathbb R$ be a derivator. Fix $x_0\in [a,b]$. We have that:
 \begin{enumerate}[itemsep=0mm]
 \item[(1)] For $x\geq x_0$ and $n\in\mathbb N$, $g_{x_0,n}(x)\geq 0$.
 \item[(2)] For $x\leq x_0$ and $n\in\mathbb N$, $g_{x_0,2n}(x)\geq 0$ and $g_{x_0,2n-1}(x)\leq 0$.
 \end{enumerate}
 Also, if $g$ is continuous on $[a,b]$,
 \[
 g_{x_0,n}(x)=g_{x_0,1}^n(x)
 \]
 for all $x\in[a,b]$. Furthermore, take $x> x_0$, if $g$ is such that $g^C=0$ and $|[x_0,x)\cap D_g|=m\in \mathbb N$. Then
 \[
 g_{x_0,n}(y)=0
 \]
 for $n\geq m+1$ and $y\in[x_0,x]$.
\end{pro}

We will use the following proposition.
\begin{pro}[{\cite[Proposition 3.2]{Frigon2017}}]
	\label{continua}
	If $f:[a,b] \to \mathbb{F}$ is $g$--continuous on $[a,b]$, then
	\begin{enumerate}[noitemsep, itemsep=.1cm]
		\item[~1.] $f$ left-continuous at every point $x\in [a,b]$;
		\item[~2.] if $g$ is continuous at $x\in [a,b]$, then $f$ is continuous at $x\in [a,b]$;
		\item[~3.] if $g$ is constant on some $[\alpha,\beta]$, then $f$ is constant on $[\alpha,\beta]\cap [a,b]$.
	\end{enumerate}
\end{pro}

In particular all $g^C$--continuous functions are continuous. The fact that $f$ must be constant where $g$ is constant allows to construct a bijection between
$g$--continuous functions on $[a,b]$ and usually continuous functions on $g([a,b])$.

\begin{lem}
 \label{compo}
 Every (uniformly) $g$--continuous function is of the form $f\circ g$ where $f$ is a (uniformly) continuous function defined on $g([a,b])$.
\end{lem}
\begin{proof}
 Clearly, $f\circ g$ is (uniformly) $g$--continuous if $f$ is a (uniformly) continuous function defined on $g([a,b])$. Take any (uniformly) $g$--continuous function $h:[a,b]\to \mathbb F$.
 Define, for $y\in g([a,b])$,
 \[
 f(y): = h(x),\quad\text{if $g(x)=y$.}
 \]
 $f$ is well-defined. Indeed, let $x'\in [a,b]$ be such that $g(x')=g(x)$, then $g$ is constant on the line segment that connects $x$ and $x'$ and, hence, by Proposition~\ref{continua},
 $h(x)=h(x')$. The above means $ f \circ g = h$.

 Furthermore, let us show that $f$ is continuous. Fix $\varepsilon>0$, and take $y=g(x)$ for $x\in[a,b]$. There exists $\delta>0$ such that
 \[
 x'\in B_g(x,\delta) \Rightarrow |h(x)- h(x')|<\varepsilon.
 \]
 Take $y'\in g([a,b])\cap (y-\delta,y+\delta)$. Let $x'\in [a,b]$ such that $g(x')=y'$. Clearly, $|g(x)-g(x')|=|y-y'|<\delta$, so $x'\in B_g(x,\delta)$ and
 \[
 |f(y)-f(y')|=| f\circ g(x)- f\circ g(x')|=|h(x)-h(x') |<\varepsilon.
 \]
 Thus, $f$ is continuous. The uniformly continuous case is analogous.
 \end{proof}
\begin{rem}
 A similar result for the uniformly continuous case with $f$ defined on the smallest interval containing $g([a,b])$ appears in \cite[Theorem 3.13]{VillaTojo}. Here we give an elementary proof.
\end{rem}
\begin{dfn}
 Define $\operatorname*{UC}_g([a,b])$ to be the Banach space of uniformly $g$--continuous functions normed with the supremum norm. Also, for the identity derivator $g=\operatorname*{id}$, we denote
 \[
 \operatorname{UC}_{\operatorname*{id}}([a,b])=\operatorname*{UC}([a,b]).
 \]
\end{dfn}

\begin{cor}
 \label{isometricgimagen}
 The map
 \begin{equation*}
 \label{isomorfismoUC}
 \begin{tikzcd}[row sep= 0em]
 \operatorname*{UC}(g([a,b])) \arrow{r} & \operatorname*{UC}_g([a,b]) \\
 f \arrow[mapsto]{r} & f\circ g
 \end{tikzcd}
 \end{equation*}
 is an isometric isomorphism.
\end{cor}

\begin{proof}
 It is clearly linear and isometric. By Lemma~\ref{compo} it is bijective.
\end{proof}


From now on, whenever we refer to the ``Weierstrass Approximation Theorem'' or simply ``the Weierstrass Theorem'', we mean the following famous result.

\begin{thm}[{\cite{1885}}]
 For every $\varepsilon<0$ and $f\in\operatorname{C}([a,b])$, there exists a polynomial $p$ such that
 \[
 \operatorname*{sup}_{x\in[a,b]}|f(x)-p(x)|\leq \varepsilon.
 \]
\end{thm}

We will sometimes work on the Hilbert space of square $g$--integrable functions.

\begin{dfn}
 We denote by $L^2_g([a,b])$ the Hilbert space of $g$--measurable functions $f:[a,b]\to \mathbb F$ that satisfy
 \[
 \int_{[a,b)} |f|^2\operatorname{d}\mu_g <\infty.
 \]
 We will denote the norm of this Hilbert space by
 \[
 \|f\|_{L^2_g} = \left(\int_{[a,b)} |f|^2\operatorname{d}\mu_g \right)^{\frac{1}{2}}.
 \]
\end{dfn}

The following lemma is needed.

\begin{lem}
	\label{lemmaintervalos}
	Let $x_0\in [a,b]$. Consider the interval
	\[I=g^{-1}(\{ g(x_0) \}) = \{ x\in [a,b] : g(x)=g(x_0) \}. \]
	Then, $\mu_g(I) = \Delta g(\operatorname{max} I )$.
\end{lem}
\begin{proof}
	Denote $I = g^{-1}(\{g(x_0)\})\subset[a,b]$. Clearly $I$ is an interval. Either $I$ is degenerate or non--degenerate, in fact,
	\[
	I=\begin{dcases}
		[\alpha,\beta], & I \text{ non--degenerate and } \alpha=\operatorname*{min} I, \beta=\operatorname*{max} I,\\
		(\alpha,\beta], & I \text{ non--degenerate and } \alpha=\operatorname*{inf} I, \beta=\operatorname*{max} I,\\
		\{x_0\}, & I \text{ degenerate}.
	\end{dcases}
	\]
	Let $\beta=\operatorname*{sup} I$ and $\alpha=\operatorname*{inf} I$. Given the left--continuity of $g$, $\beta$ always belongs to $I$. There are two cases.

	If $\alpha \in I$, we have that
	\begin{align*}
		\mu_g(I) = \mu_g([\alpha,\beta]) = g(\beta ^+)-g(\alpha)=g(\beta^+)-g(\beta) = \Delta g(\beta),
	\end{align*}
	since $g$ is constant on $I$.

	If $\alpha \notin I$, we have that
	\begin{align*}
		\mu_g(I) = \mu_g((\alpha,\beta]) = g(\beta ^+)-g(\alpha^+)=g(\beta^+)-g(\beta) = \Delta g(\beta),
	\end{align*}
	since $g$ is constant on $I$.

	Notice that if $\beta = b$ we cannot take the right hand side limit $g(b^+)$. However, in that case $\mu_g(I)=0$.
\end{proof}

We present now a useful density result.

\begin{pro}
	\label{density}
	The space $\operatorname*{UC}_g([a,b])$ is dense on $L^2_g([a,b])$. Furthermore, let $x_0\in [a,b]$ be such that the set
	\[
	g^{-1}(\{ g(x_0) \}) = \{ x\in [a,b] : g(x)=g(x_0) \}
	\]
	has null $g$--measure. Then the subspace
	\[
	\{ f\in \operatorname{UC}_g([a,b]) \,: \,f(x_0)= 0  \}
	\]
	is dense on $L^2_g([a,b])$.
\end{pro}
\begin{proof}
	We will show first that $\operatorname*{UC}_g([a,b])$ is dense on $L^2_g([a,b])$. Denote $
	C = \overline{\operatorname{UC}_g([a,b])}.$
	We only have to prove that the indicator functions $1_{[x_1,x_2)}$ belong to $C$ for all $a\leq x_1< x_2\leq b$ such that $g(x_2)-g(x_1)>0$ (so $\mu_g([x_1,x_2))>0$). Then, for any $g$--measurable set $E\subset [a,b]$, $1_E$ belongs to $C$ as it can be approximated with functions of the type of $1_{[x_1,x_2)}$. Finally, any $f\in L^2_g([a,b])$ is the limit of simple functions and the result follows.

	Fix $x_1,x_2\in[a,b]$ such that $g(x_2)-g(x_1)>0$. Let $n\in \mathbb N$. Define $f:\mathbb R \to [0,1]$ as follows:
	\begin{equation}
		\label{funcionmesetadef}
	f(y) = \begin{cases}
		1, & g(x_1)\leq y\leq g(x_2)-\frac{1}{n},\\
		0, & y\leq g(x_1)-\frac{1}{n} \text{ or }y\geq g(x_2),\\
		n(y-(g(x_1)-\frac{1}{n})), & g(x_1)-\frac{1}{n}<y<g(x_1),\\
		-n(y-g(x_2))), & g(x_2)-\frac{1}{n}<y<g(x_2),
	\end{cases}
	\end{equation}

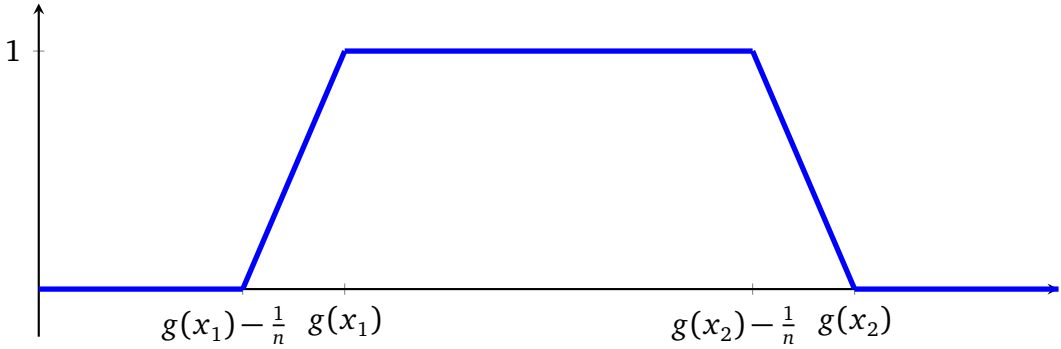
\begin{figure}[h]
	\centering
\begin{tikzpicture}
\begin{axis}[
    axis lines=middle,
    ymin=-0.2, ymax=1.2,
    xmin=0, xmax=10,
    samples=200,
    domain=0:10,
    thick,
    width=15cm,
    height=6cm,
    xtick={2,3,7,8},
    xticklabels={\hspace{-1.2em}$g(x_1)-\frac{1}{n}$, $g(x_1)$, \hspace{-1.2em}$g(x_2)-\frac{1}{n}$, $g(x_2)$},
    ytick={0,1},
]
\addplot[blue,line width=2pt, domain=0:2]{0};
\addplot[blue,line width=2pt,domain=2:3]{x - 2};
\addplot[blue,line width=2pt,domain=3:7]{1};
\addplot[blue,line width=2pt,domain=7:8]{-x + 8};
\addplot[blue,line width=2pt,domain=8:10]{0};

\end{axis}
\end{tikzpicture}
\caption{Example of function $f$ in \eqref{funcionmesetadef}.}
\label{funcionmeseta}
\end{figure}
\noindent for $y\in\mathbb R$ ---see Figure \ref{funcionmeseta}. We have that $f$ is uniformly continuous on $\mathbb R$. Therefore, by Lemma~\ref{compo}, $f\circ g$ is uniformly $g$--continuous on $[a,b]$. Let $h(x): = |f\circ g- 1_{[x_1,x_2)}|(x)$, for $x\in[a,b]$. Note that $\|h\|_\infty\leq 2$. We are going to show that, for large enough $n$, $\|h\|_{L^2_g}$ is as small as we want. Consider the sets
	\[\begin{array}{l}
		S_1 = \{z\in [a,b]: g(z)\leq g(x_1)-\frac{1}{n} \text{ or }g(z)\geq  g(x_2)\},\\
		S_2 = \{z\in [a,b]: g(x_1) \leq g(z)\leq g(x_2) - \frac{1}{n}\}\cap [x_1,x_2).
	\end{array}
	\]
	Let us prove that $h(x)=0$ for $g$--a.a. $x\in S_1\cup S_2$.

\textbullet	If $x\in S_1$, we have two cases: either $g(x)\leq g(x_1)-\frac{1}{n}$ and, therefore, $x\notin [x_1,x_2)$ and $h(x)=0$; or $g(x)\geq  g(x_2)$, in which case, we have two cases:

	\noindent (a) $g(x)> g(x_2)$, in which case $x\notin [x_1,x_2)$ and $h(x)=0$, or

	\noindent (b) $g(x)= g(x_2)$. If $x\notin [x_1,x_2)$ we have that $h(x)=0$. However, $x$ may belong to $[x_1,x_2)$. In that case, following Lemma~\ref{lemmaintervalos}, $\mu_g(g^{-1}(\{g(x_2)\})\cap [x_1,x_2))=0$.

\textbullet	If $x\in S_2$, then $h(x)=0$ and the claim follows.

It is only left to study the set $[a,b]\setminus (S_1\cup S_2)$. Denote
	\[
	\begin{array}{l}
		A_1 = \{z\in [a,b]:  g(x_1)-\frac{1}{n}<g(z)<g(x_1)\},\\
		A_2 = \{z\in [a,b]:  g(x_2)-\frac{1}{n}<g(z)<g(x_2)\},\\
		A_3 = \{z\in [a,b]:  g(x_1) \leq g(z)\leq g(x_2) - \frac{1}{n}\}\setminus [x_1,x_2).
	\end{array}
	\]
	Hence, $[a,b]\setminus (S_1\cup S_2) = A_1\cup A_2\cup A_3$. Note also that
	\[
	A_3 = g^{-1}(\{g(x_1)\})\cap [a,x_1).
	\]
	Thus, $A_3$ has null $g$--measure. We show $\mu_g(A_1)\leq 1/n$. If $A_1$ is empty the claim follows. If not denote $s = \operatorname*{inf} A_1$. Consider again two cases:

	\textbullet	$s\in A$. Then, $A_1\subset [s,x_1)$ and hence $\mu_g(A_1)\leq g(x_1)-g(s)<1/n$.

	\textbullet	$s\notin A$. Then, $A_1\subset (s,x_1)$ and hence $\mu_g(A_1)\leq g(x_1)-g(s^+)<1/n$.

	Similarly, we can show that $\mu_g(A_2)\leq 1/n$. Whence, we have that $\mu_g([a,b]\setminus (S_1\cup S_2))\leq 2/n$. Thus,
	\begin{align*}
		\|f\circ g - 1_{[x_1,x_2)}\|^2_{L^2_g} &=\|h\|^2_{L^2_g} = \int_{[a,b)} h^2\operatorname*{d}\mu_g \\
		& =\int_{[a,b)\setminus (S_1\cup S_2)} h^2\operatorname*{d}\mu_g\leq \|h\|^2_\infty\frac{2}{n}\leq \frac{8}{n}\xrightarrow[n\to \infty]{}0.
	\end{align*}

	Let us show that $\{ f\in \operatorname{UC}_g([a,b]) \,: \,f(x_0)= 0  \}$ is dense on $\operatorname*{UC}_g([a,b])$ in $L^2_g([a,b])$ if $\mu_g(g^{-1}(\{ g(x_0) \}))=0$. Let $f\in\operatorname*{UC}_g([a,b])$ and denote $I=g^{-1}(\{g(x_0)\})$, $\alpha = \operatorname*{inf} I$ and $\beta = \operatorname*{max} I$. We have that $\alpha \in I \Leftrightarrow \Delta g(\alpha)=0$, otherwise it would contradict Lemma~\ref{lemmaintervalos}.

Henceforth we impose the conditions $a\ne\a\ne\b\ne b$. These assumptions streamline the exposition; if any of them fails, the construction carries over after straightforward adjustments.

	Let $n\in \mathbb N$ be such that $n>\max\{\frac{1}{\a-a},\frac{1}{b-\b}\}$. Note that, given the definition of $\beta$, $g(\beta +1/n)\ne g(\beta)=g(x_0)$. The same for $\alpha$ when $\alpha\notin D_g$. With these considerations in mind, if $\alpha\notin D_g$, define
	\begin{equation}
		\label{alphanodg}
		h(x) = \begin{dcases}
		f(x), & a\leq x\leq \alpha -\frac{1}{n},\\
		f\left(\alpha -\frac{1}{n}\right)\frac{g(x)-g(\alpha)}{g\left(\alpha -\frac{1}{n}\right)-g(\alpha)},& \alpha -\frac{1}{n}<x\leq \alpha,\\
		0, & \alpha < x \leq \beta,\\
		f\left(\beta +\frac{1}{n}\right)\frac{g(x)-g(\beta)}{g\left(\beta +\frac{1}{n}\right)-g(\beta)},& \beta <x\leq \beta+\frac{1}{n},\\
		f(x), & \beta+\frac{1}{n} < x\leq b.
	\end{dcases}
	\end{equation}

\begin{figure}[h]
	\centering
\begin{tikzpicture}
\begin{axis}[
    axis lines=middle,
    ymin=-1.2, ymax=1.2,
    xmin=0, xmax=10,
    samples=200,
    domain=0:10,
    thick,
    width=15cm,
    height=6cm,
	ytick={0},
	xtick={0},
	declare function={
        f(\x) = sin(deg(\x));
		gu(\x) = x-6;
		gd(\x) = x-4;
		s(\x) = -1 * (\x <= 3.5);
    },
]
\addplot[blue,line width=2pt,domain=0:3]{f(x)};
\addplot[blue,line width=2pt,domain=3:4]{(f(3)/(-2))*(gd(x)+s(x))};
\addplot[blue,line width=2pt,domain=4:6]{0};
\addplot[blue,line width=2pt,domain=6:7]{(f(5))*gu(x)};
\addplot[blue,line width=2pt,domain=7:10]{f(x-2)};
\addplot[dashed,line width=2pt,domain=3:4]{f(x)};
\addplot[dashed,line width=2pt,domain=6:7]{f(x-2)};
\addplot[dashed,line width=2pt,domain=4:6]{f(4)};
\end{axis}
\end{tikzpicture}
\caption[]{\protect{Example of function $h$ in \eqref{alphanodg} (solid line) with derivator
	
\begin{minipage}{\textwidth}
\[
g(x) = \begin{cases}
	x, & x\in[0,3.5],\\
	x+1, & x\in(3.5,4],\\
	5, & x\in[4,6],\\
	x-1, & x\in(6,10].
\end{cases}
\]
\end{minipage}
	We choose $n=1$ and $x_0=5$, so $\alpha = 4$ and $\beta = 6$. In a dashed line, the function $f$.
}}
\label{alphanofgfig}
\end{figure}
---see Figure~\ref{alphanofgfig}. If $\alpha \in D_g$, define
	\begin{equation}
		\label{alphasidg}
		h(x) = \begin{dcases}
		f(x), & a\leq x\leq \alpha,\\
		0, & \alpha < x \leq \beta,\\
		f\left(\beta +\frac{1}{n}\right)\frac{g(x)-g(\beta)}{g\left(\beta +\frac{1}{n}\right)-g(\beta)},& \beta <x\leq \beta+\frac{1}{n},\\
		f(x), & \beta+\frac{1}{n} < x\leq b.
	\end{dcases}
	\end{equation}

\begin{figure}[t]
	\centering
\begin{tikzpicture}
\begin{axis}[
    axis lines=middle,
    ymin=-1.2, ymax=1.2,
    xmin=0, xmax=10,
    samples=200,
    domain=0:10,
    thick,
    width=15cm,
    height=6cm,
    ytick={0},
	declare function={
        f(\x) = sin(deg(\x));
		gu(\x) = x-6;
		s(\x) = 1 * (\x > 6.5);
    },
]
\addplot[blue,line width=2pt,domain=0:4]{f(x)};
\addplot[blue,line width=2pt,domain=4:6]{0};
\addplot[blue,line width=2pt,domain=6:7]{(f(7)/2)*(gu(x)+s(x))};
\addplot[blue,line width=2pt,domain=7:10]{f(x)};

\addplot[dashed,line width=2pt,domain=6:7]{f(x)};
\addplot[dashed,line width=2pt,domain=4:6]{f(6)};

\end{axis}
\end{tikzpicture}
\caption[]{\protect{Example of function $h$ in \eqref{alphasidg} (solid line) with derivator
		
		\begin{minipage}{\textwidth}
			\[
			g(x) = \begin{cases}
				x, & x\in[0,4],\\
				5, & x\in(4,6],\\
				x-1, & x\in[6,6.5],\\
				x, & x\in(6.5,10].
			\end{cases}
			\]
		\end{minipage}
		We choose $n=1$ and $x_0=5$, so $\alpha = 4$ and $\beta = 6$. In a dashed line, the function $f$.
}}\label{alphasidgfig}
\end{figure}
	---see Figure~\ref{alphasidgfig}. We now want to show that $h\in\{ f\in \operatorname{UC}_g([a,b]) \,: \,f(x_0)= 0  \}$ and $\|h-f\|_{L^2_g}$ is as small as we want for $n$ large enough.

	Clearly, $h$ is $g$--continuous for $x\in[a,b]\bs\{\alpha-1/n,\alpha,\beta,\beta+1/n\}$. We show that $h$ is $g$--continuous on $\beta+1/n$ and in $\alpha$ when $\alpha\in D_g$. The rest of the cases are analogous. Denote $t = \beta+1/n$ and let $\varepsilon>0$. There exists $\delta'>0$ such that
	\[
	|f(x)-f(t)|<\varepsilon \text{ for all }x\in [a,b] \text{ such that }|g(x)-g(t)|<\delta'.
	\]
 If $|f(t)|=0$ let $\delta =\delta'$, otherwise, take  $\delta\in(0, \min\{\delta',\varepsilon|f(t)|^{-1}|g(t)-g(\beta)|,|g(t)-g(\beta)|\})$. Note that, $B_g(t,\delta)\subset (\beta,b]$. Let $x\in[a,b]$ be such that $|g(x)-g(t)|<\delta$. If $x\leq t$,
	\[
	|h(x)-h(t)| = \left|f(t)\frac{g(x)-g(t)}{g(t)-g(\beta)} \right|< \varepsilon.
	\]
	If $x>t$,
	\[
	|h(x)-h(t)| = |f(x)-f(t)|< \varepsilon.
	\]
	Thus, $h$ is $g$--continuous at $t$.

	Now, take $\delta\in(0,\Delta g(\alpha))$. If $x\in[a,b]$ is such that $|g(x)-g(\alpha)|<\delta$ then $x$ must belong to $[a,\alpha]$. But  $h$ and $f$ are equal on this interval so $h$ is $g$--continuous at $\alpha$.

	Since $h$ is $g$--continuous on all $[a,b]$ and regulated, it is uniformly $g$--continuous, see \cite[Corollary~3.8]{VillaTojo}. Furthermore, $h(x_0)=0$. If $\alpha\notin D_g$,
	\begin{align*}
		\|h-f\|^2_{L^2_g} &= \int_{[a,b)} |h-f|^2 \operatorname*{d}\mu_g = \int_{(\alpha-\frac{1}{n},\beta+\frac{1}{n})} |h-f|^2 \operatorname*{d}\mu_g \leq 4\|f\|^2_\infty \left(g\left(\beta+\frac{1}{n}\right)-g\left(\alpha-\frac{1}{n}\right)\right).
	\end{align*}
Recall we assumed $g^{-1}(\{g(x_0)\})$ has null $g$--measure and this means $\Delta g(\beta)=0$ by Lemma~\ref{lemmaintervalos}.	Now, $\lim_{n\to\infty} g(\alpha-1/n)=g(\alpha)$ by the left--continuity of $g$, $g(\alpha)=g(\beta)$ since $\alpha \in I$ ($\alpha\notin D_g$) and $\lim_{n\to\infty} g(\beta+1/n)=g(\beta)$ since $\Delta g(\beta)=0$.  Therefore, $\|h-f\|_{L^2_g}$ can be made as small as we want for $n$ sufficiently large.

	If $\alpha\in D_g$,
	\begin{align*}
		\|h-f\|^2_{L^2_g} &= \int_{[a,b)} |h-f|^2 \operatorname*{d}\mu_g = \int_{(\alpha,\beta+\frac{1}{n})} |h-f|^2 \operatorname*{d}\mu_g\leq 4\|f\|^2_\infty \left(g\left(\beta+\frac{1}{n}\right)-g(\alpha^+)\right),
	\end{align*}
	and $\lim_{n\to\infty}g(\beta+\frac{1}{n})=g(\beta)=g(\alpha^+)$. Again, as $n$ tends to infinity, $\|h-f\|_{L^2_g}$ tends to zero. Thus, $\{ f\in \operatorname{UC}_g([a,b]) \,: \,f(x_0)= 0  \}$ is dense on $L^2_g([a,b])$.
\end{proof}

\section{Routes of proof}
\label{seccionrutas}
In this section, we discuss the proofs of the Weierstrass approximation theorem in classical analysis and how they adapt to the Stieltjes setting. We show first a trivial case where one obtains density via straightforward application of the classical theorem---cf. \cite[Theorem 3.17]{VillaTojo}.

\begin{pro}
 \label{casotontogcontinuo}
 Let $g:[a,b]\to \mathbb R$ be a continuous derivator. Then, ${\rm P}_g$ is dense in $\operatorname{UC}_{g}([a,b])$.
\end{pro}
\begin{proof}
 By Proposition~\ref{propiedadesmonomios}, we know that $g_{a,n}(x) = (g(x) - g(a))^n$ for all $x \in \mathbb{R}$ and $n \in \mathbb{N}$. Therefore, any $g$--polynomial can be written as $p \circ g$, where $p$ is an ordinary polynomial, and conversely, any composition $p \circ g$ is a $g$--polynomial. According to Corollary~\ref{isometricgimagen}, the space $\operatorname{UC}(g([a,b])) = \operatorname{UC}([g(a), g(b)])$ is isometric to $\operatorname{UC}_g([a,b])$ via composition with $g$. The Weierstrass Approximation Theorem ensures that ordinary polynomials are dense in $\operatorname{UC}([g(a), g(b)])$, which implies that ${\rm P}_g$ is dense in $\operatorname{UC}_g([a,b])$.
\end{proof}

The challenge arises when $g$ is not continuous, in which case ${\rm P}_g$ no longer is an algebra with the product of functions: the product of two $g$--polynomials is, in general, not a $g$--polynomial. This not only prevents us from applying the Stone-Weierstrass Theorem, but also invalidates many traditional strategies used in polynomial approximation.

\begin{thm}[{\cite{Stone1948}}]
 Suppose $X$ is a compact Hausdorff space and $A$ is a subalgebra of $\operatorname{C}(X, \mathbb{R})$ that contains a non-zero constant function. Then $A$ is dense in $\operatorname{C}(X, \mathbb{R})$ if and only if it separates points.
\end{thm}

Another difficulty in adapting Weierstrass's theorem to Stieltjes calculus is that classical approximation techniques do not readily extend to derivators. Weierstrass's original proof relies on convolution with a Gaussian heat kernel—essentially, he uses analytic functions that form an approximate identity. Properly truncating the power series of the Gaussian kernel and interchanging it with convolution yields a sequence of polynomials that converge uniformly (see \cite{1885}). Sadly, convolution has not yet been developed within the framework of Stieltjes calculus.

An alternative approach is based on Bernstein polynomials \cite{Bernstein}. However, the binomial theorem in Stieltjes calculus \cite[Proposition 3.14]{Cora2023} is insufficient for defining Bernstein polynomials in the same way as in the classical setting. Moreover, transferring the coefficients of the Bernstein polynomials directly onto the corresponding $g$--polynomials loses many essential properties required for the proof.

Both Bernstein's and analytic approximations to the identity type proofs can be viewed as sequences of positive linear operators on the space of continuous functions. Bohman and Korovkin independently established the following fundamental result \cite{Korovkin1953, Bohman1952}.

\begin{thm}
 Let $(L_n)_{n \in \mathbb{N}}$ be a sequence of positive linear operators defined on $\operatorname{C}([a,b])$. If $L_n(f) \to f$ uniformly on $[a,b]$ for $f \in \{1, x, x^2\}$, then $L_n(f) \to f$ uniformly for all $f \in \operatorname{C}([a,b])$.
\end{thm}

Both Weierstrass's and Bernstein's proofs fall under this theorem. Examples of such operator sequences include Fejér kernels (used for trigonometric polynomial density) and those used in Landau's approach of Weierstrass's Theorem \cite{Landau1908}. Even if an analogue of Korovkin’s theorem were available for Stieltjes calculus, one would still need to construct an appropriate sequence of operators.

A promising strategy is to study the closure of ${\rm P}_g$ in the Hilbert space $L^2_g([a,b])$, where the inner product simplifies distance calculations. Let us first show the equivalence of the two density problems. The following lemma is useful.

\begin{lem}
 \label{lintrozos}
 Denote the vector space
 \[ M = \{ f\in \operatorname{AC}_g([a,b]) : f' \in L^2_g([a,b])\}\subset\operatorname{UC}_g([a,b]).\]
 Then, $M$ is dense in $\operatorname{UC}_g([a,b])$.
\end{lem}

\begin{proof}
 Let $f\in \operatorname{UC}_{g}([a,b])$ and let $\varepsilon>0$. By uniform $g$--continuity, there exists $\delta>0$ such that
 \[
 |g(x)-g(y)|\leq \delta \text{ implies } |f(x)-f(y)|\leq \varepsilon
 \]
 for all $x,y\in[a,b]$. Let $a = x_0<\dots<x_m=b$ be a partition of $[a,b]$ such that
 \begin{equation}
	\label{ec:part}
	\operatorname*{sup}_{i=1,\dots,m}|g(x_i)-g(x_{i-1}^+)|\leq \delta.
 \end{equation} 
 We will later prove the existence of such a partition. Uniformly $g$--continuous functions are regulated \cite[Corollary 3.8]{VillaTojo}. Consider the following function in $L^2_g([a,b])$,
 \begin{equation}
	\label{derivadal2fun}
	 h(x)=\begin{dcases}
 \frac{f(x_i)-f(x_{i-1}^+)}{g(x_i)-g(x_{i-1}^+)}, & \text{ if } x\in (x_{i-1},x_i) \text{ and } g(x_i)-g(x_{i-1}^+)>0,\\
 f'_g(x_i), & \text{ if } x = x_i \text{ and } \Delta g(x)>0,\\[ .4em]
 0, & \text{ otherwise}.
 \end{dcases}
 \end{equation}
---see Figure~\ref{derivadal2funfig} for an example.
\begin{figure}[h]
\centering
\begin{tikzpicture}
\begin{axis}[
    axis lines=middle,
    ymin=-2.1, ymax=2.1,
    xmin=0, xmax=10,
    samples=200,
    domain=0:10,
    thick,
    width=15cm,
    height=6cm,
    ytick={0},
	xtick={0},
	declare function={
        f(\x) = sin(deg(\x));
    },
]
\addplot[black,line width=2pt,domain=0:1]{f(1)-f(0)};
\addplot[black,line width=2pt,domain=1:2]{f(2)-f(1)};
\addplot[black,line width=2pt,domain=2:3]{f(3)-f(2)};
\addplot[black,line width=2pt,domain=3:4]{(f(4)-1)-(f(3)-1)};
\addplot[black,line width=2pt,domain=4:5.5]{(f(5)-1)-(f(4)-1)};
\addplot[black,line width=2pt,domain=5.5:6]{2*((f(5.5)-1)-((f(5)-1)))};
\addplot[black,line width=2pt,domain=6:7]{(f(6.5)+1)-(f(5.5)+1)};
\addplot[black,line width=2pt,domain=7:8]{(f(7.5)+1)-(f(6.5)+1)};
\addplot[black,line width=2pt,domain=8:9]{(f(8.5)+1)-(f(7.5)+1)};
\addplot[black,line width=2pt,domain=9:10]{(f(9.5)+1)-(f(8.5)+1)};

\addplot[black, only marks, mark=*, mark size=1.5pt] coordinates {(3,-1)};
\addplot[black, only marks, mark=*, mark size=1.5pt] coordinates {(6,2)};
\addplot[dashed,line width=2pt,domain=0:3]{f(x)};
\addplot[dashed,line width=2pt,domain=3:4.5]{f(x)-1};
\addplot[dashed,line width=2pt,domain=4.5:5]{f(4.5)-1};
\addplot[dashed,line width=2pt,domain=5:6]{f(x-0.5)-1};
\addplot[dashed,line width=2pt,domain=6:10]{f(x-0.5)+1};
\end{axis}
\end{tikzpicture}
\caption[]{\protect{Example of function $h$ in \eqref{derivadal2fun} (solid line) with derivator $g=g^*+g^{**}$, where
				
			\begin{minipage}{\linewidth}
				\begin{equation*}
					g^*(x) = \begin{dcases}
						x, & x\in[0,4.5],\\
						4.5, & x\in(4.5,5],\\
						x-1.5, & x\in(5,10],
					\end{dcases}
					\quad
					g^{**}(x) = \begin{dcases}
						0, & x\in[0,3],\\
						1, & x\in(3,6],\\
						2, & x\in(6,10].
					\end{dcases}
				\end{equation*}
			\end{minipage}
			We choose as a partition $\{0,1,2,\dots,10\}$, which satisfies~\eqref{ec:part} with $\delta = 1$. Plotted with a dashed line, the function $f$.}}
\label{derivadal2funfig}
\end{figure}

Note that, for all $i\in\{1,\dots,m\}$, $g(x_i)-g(x_{i-1}^+)=0$ if and only if $(x_{i-1},x_i)$ has null $g$--measure. Define the function
 \begin{equation}
	\label{integralderfunl2}
	L(x) = f(x_0)+\int_{[x_0,x)}h\operatorname{d}\mu_g.
 \end{equation}
 We have that $L\in M$ and moreover, it can be written explicitly as
 \[
 L(x) = \begin{dcases}
 f(x_0), & \text{ if } x=x_0,\\[ .7em]
 f(x_{i-1}^+), & \text{ if } x\in (x_{i-1},x_i]\text{ and } g(x_i)-g(x_{i-1}^+)=0,\\
 f(x_{i-1}^+)+(g(x)-g(x_{i-1}^+))\frac{f(x_i)-f(x_{i-1}^+)}{g(x_i)-g(x_{i-1}^+)}, & \text{ if } x\in (x_{i-1},x_i]\text{ and } g(x_i)-g(x_{i-1}^+)>0.
 \end{dcases}
 \]
---see Figure~\ref{integralderfunl2fig}.

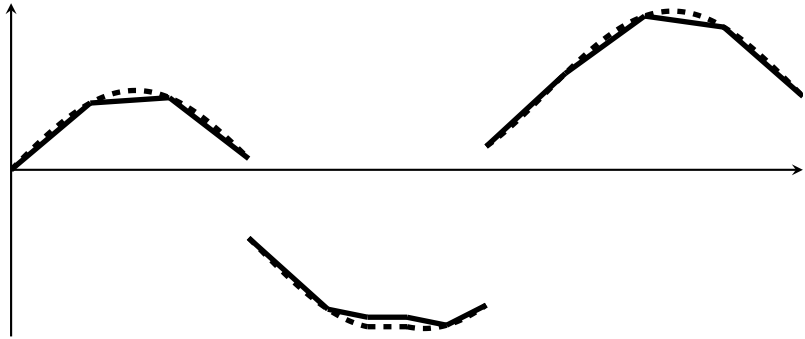
\begin{figure}[h]
\centering
\begin{tikzpicture}
\begin{axis}[
    axis lines=middle,
    ymin=-2.1, ymax=2.1,
    xmin=0, xmax=10,
    samples=200,
    domain=0:10,
    thick,
    width=12cm,
    height=6cm,
    ytick={0},
	xtick={0},
	declare function={
        f(\x) = sin(deg(\x));
    },
]
\addplot[black,line width=2pt,domain=0:1]{(f(1)-f(0))*x+f(0)};
\addplot[black,line width=2pt,domain=1:2]{(f(2)-f(1))*(x-1)+f(1)};
\addplot[black,line width=2pt,domain=2:3]{(f(3)-f(2))*(x-2)+f(2)};
\addplot[black,line width=2pt,domain=3:4]{((f(4)-1)-(f(3)-1))*(x-3)+f(3)-1};
\addplot[black,line width=2pt,domain=4:4.5]{((f(5)-1)-(f(4)-1))*(x-4)+(f(4)-1)};
\addplot[black,line width=2pt,domain=4.5:5]{((f(5)-1)-(f(4)-1))*(4.5-4)+(f(4)-1)};
\addplot[black,line width=2pt,domain=5:5.5]{((f(5)-1)-(f(4)-1))*(x-0.5-4)+(f(4)-1)};
\addplot[black,line width=2pt,domain=5.5:6]{(2*((f(5.5)-1)-((f(5)-1))))*(x-5.5)+(f(5)-1)};
\addplot[black,line width=2pt,domain=6:7]{((f(6.5)+1)-(f(5.5)+1))*(x-6)+(f(5.5)+1)};
\addplot[black,line width=2pt,domain=7:8]{((f(7.5)+1)-(f(6.5)+1))*(x-7)+(f(6.5)+1)};
\addplot[black,line width=2pt,domain=8:9]{((f(8.5)+1)-(f(7.5)+1))*(x-8)+(f(7.5)+1)};
\addplot[black,line width=2pt,domain=9:10]{((f(9.5)+1)-(f(8.5)+1))*(x-9)+(f(8.5)+1)};

\addplot[dashed,line width=2pt,domain=0:3]{f(x)};
\addplot[dashed,line width=2pt,domain=3:4.5]{f(x)-1};
\addplot[dashed,line width=2pt,domain=4.5:5]{f(4.5)-1};
\addplot[dashed,line width=2pt,domain=5:6]{f(x-0.5)-1};
\addplot[dashed,line width=2pt,domain=6:10]{f(x-0.5)+1};
\end{axis}
\end{tikzpicture}
\caption{Example of function $L$ in \eqref{integralderfunl2} (solid line) where we integrate the example function of Figure~\ref{derivadal2funfig}. In a dashed line, again, the function $f$.}
\label{integralderfunl2fig}
\end{figure}

 We now show that $\|L-f\|_{\infty}\leq \varepsilon$. Let $x\in[a,b]$. If $x=x_i$ for some $i$ then $L(x)=f(x)$. If $x\in (x_{i-1},x_i)$ for some $i\in\{1,\dots,m\}$, we distinguish two cases. If $g(x_i)-g(x_{i-1}^+)=0$, then $f$ is constant and equal to $f(x_{i-1}^+)$ on $(x_{i-1},x_i)$ and hence $|L(x)-f(x)|=0$. If $(x_{i-1},x_i)$ does not have zero $g$--measure, denote $\lambda = (g(x)-g(x_{i-1}^+))(g(x_i)-g(x_{i-1}^+))^{-1}\in[0,1]$, then,
 \begin{align*}
 |L(x)-f(x)| &= |(1-\lambda)f(x_{i-1}^+)+ \lambda f(x_{i})-f(x)|\\
 &\leq (1-\lambda)|f(x_{i-1}^+)-f(x)|+ \lambda |f(x_{i})-f(x)|\leq \varepsilon,
 \end{align*}
 since $|g(x)-g(y)|\leq \delta$ for all $x,y\in (x_{i-1},x_i]$.

It remains to prove that such a partition exists. Define $x_0=a$ and consider the interval
\[
 \{z\in[a,b]\,:\, g(z)\leq g(x_{0}^+)+\delta\}.
\]
The above set is a compact subset of $\mathbb R$ on the usual topology and hence it has a maximum. Thus, define $x_1$ to be that maximum. By construction, it holds that
 \[ g(x_{1})-g(x_{0}^+)\leq \delta. \]
 Moreover, either $x_{1} = b$ or it also holds that
 \[
 g(x_{1}^+)-g(x_{0}^+)\geq \delta.
 \]
We will construct the partition inductively. Suppose the points $x_k$ have been defined up to an index $i$ in such a way that
 \[
 g(x_{k+1})-g(x_{k}^+)\leq \delta \text{ and } g(x_{k+1}^+)-g(x_{k}^+)\geq \delta,
 \]
 for all $k\in\{0,\dots,i-1\}$. If $x_{i}=b$, we are done. Otherwise, define $x_{i+1}$ as the maximum of the set
 \[
 \{z\in[a,b]\,:\, g(z)\leq g(x_{i}^+)+\delta\}.
 \]
 By the arguments made for $x_1$ either $x_{i+1}=b$ or 
 \[
 g(x_{k+1})-g(x_{k}^+)\leq \delta \text{ and } g(x_{k+1}^+)-g(x_{k}^+)\geq \delta,
 \]
 for all $k\in\{0,\dots,i\}$. In both cases we have that $x_{i+1}$ is strictly greater than $x_i$; they cannot be the same point. Besides,
 \[
 g(x_{i+1}^+)-g(x_{0}^+) = \sum_{k=0}^{i} g(x_{k+1}^+)-g(x_{k}^+)\geq (i+1)\delta,
 \]
 which implies that $\mu_g([x_0,x_{i+1}])\geq (i+1)\delta$. Therefore, we necessarily reach $x_m=b$ in a finite number of steps. \qedhere

\end{proof}

The proof of Lemma~\ref{lintrozos} shows that every uniformly $g$--continuous function can be approximated by functions that are ``piecewise $g$--linear''. Lebesgue used a similar idea to prove Weierstrass's Theorem \cite{Lebesgue1898}, along with the fact that $|x|$ can be approximated by polynomials. Unfortunately, this argument does not carry over easily to Stieltjes calculus.

\begin{pro}
 \label{densidadequiv}
 ${\rm P}_g$ is dense in $L^2_g([a,b])$ if and only if it is dense in $\operatorname{UC}_g([a,b])$.
\end{pro}
\begin{proof}
 Assume that ${\rm P}_g$ is dense in $L^2_g([a,b])$. Let $f\in \operatorname{AC}_{g}([a,b])\subset \operatorname{UC}_{g}([a,b])$ such that $f'_g\in L^2_g([a,b])$, ($f\in M$ with the same $M$ as in Lemma~\ref{lintrozos}). Define
 \[ P(x) = f(a) + \int_{[a,x)} p_g\operatorname{d}\mu_g, \]
 where $p_g$ is some $g$--polynomial. Note that $P$ is also a $g$--polynomial. For $x\in[a,b]$, we have
 \begin{align*}
 |f(x)- P(x)| = & \left| \int_{[a,x)}f'_g\operatorname{d}\mu_g- \int_{[a,x)}p_g\operatorname{d}\mu_g\right| \leq \int_{[a,x)}|f'_g-p_g|\operatorname{d}\mu_g\\
 & \leq \int_{[a,b)}|f'_g-p_g|\operatorname{d}\mu_g \leq \|1\|_{L^2_g}\|f'_g-p_g\|_{L^2_g}.
 \end{align*}
 Since we can make $\|f'_g-p_g\|_{L^2_g}$ as small as we want, we have that $M$ is contained on the closure of ${\rm P}_g$ in $\operatorname{UC}_{g}([a,b])$. Therefore, ${\rm P}_g$ is dense in $\operatorname{UC}_{g}([a,b])$ by Lemma~\ref{lintrozos}.

 Assume now ${\rm P}_g$ is dense in $\operatorname{UC}_{g}([a,b])$. Since the $L^2_g$ norm is weaker than the supremum norm and $\operatorname{UC}_{g}([a,b])$ is dense on $L^2_g([a,b])$, recall Proposition~\ref{density}, the equivalence follows.
\end{proof}

This trick, bounding the supremum norm via the $L^2$ norm on the derivatives using absolute continuity and Holder's inequality, is often employed in the proof of Müntz's theorem \cite{Muntz1914}. This theorem characterizes when the functions $\{x^{p_k}\}_{k \in \mathbb{N}}$ are dense in $\operatorname{C}([a,b])$ for a given sequence of exponents $\{p_k\}_{k \in \mathbb{N}}$---see \cite[Theorem 11.3.4]{Davis1974} for a modern treatment. They also use inner product spaces to calculate distances explicitly.

\begin{dfn}
 Let $V$ be an inner product space and $x_1, \dots, x_n \in V$. The \emph{Gram matrix} is defined by
 \[
 \operatorname{Gram}(x_1, \dots, x_n) = \begin{pmatrix}
 \langle x_1, x_1 \rangle & \langle x_1, x_2 \rangle & \cdots & \langle x_1, x_n \rangle \\
 \langle x_2, x_1 \rangle & \langle x_2, x_2 \rangle & \cdots & \langle x_2, x_n \rangle \\
 \vdots & \vdots & \ddots & \vdots \\
 \langle x_n, x_1 \rangle & \langle x_n, x_2 \rangle & \cdots & \langle x_n, x_n \rangle
 \end{pmatrix}.
 \]
 We denote its determinant by $\operatorname{gram}(x_1, \dots, x_n)$.
\end{dfn}

This determinant plays an important role in approximation theory, as shown in the following result.

\begin{lem}[{\cite[Theorem 8.7.4]{Davis1974}}]
 \label{distanciaprodescalar}
 Let $V$ be an inner product space and $x_1, \dots, x_n \in V$ be linearly independent. Denote $X = \operatorname{span}\{x_1, \dots, x_n\}$. Then, for any $v \in V$,
 \[
 d(v, X)^2 = \frac{\operatorname{gram}(v, x_1, \dots, x_n)}{\operatorname{gram}(x_1, \dots, x_n)},
 \]
 where $d(v, X)$ denotes the distance from $v$ to the subspace $X$.
\end{lem}

 Proposition~\ref{densidadequiv} gives us the following characterization.

\begin{thm}
 \label{equivalenciagpolinomiosdensidad}
 ${\rm P}_g$ is dense in $\operatorname{UC}_g([a,b])$ if and only if the indicator functions
 \[
 1_{\{x\}}(t) = \begin{dcases}
 1, & \text{if } t = x, \\
 0, & \text{if } t \neq x,
 \end{dcases}
 \]
 for $x \in D_g$ belong to the closure of ${\rm P}_g$ in $L^2_g([a,b])$.
\end{thm}
\begin{proof}
 By Proposition~\ref{densidadequiv}, the necessity is evident. Let us now prove the sufficiency. The vector space generated by $\{g^n\}_{n=0}^{\infty}$ is dense in $\operatorname{UC}_g([a,b])$, see \cite[Theorem 3.17]{VillaTojo}. We will show that we can uniformly approximate these functions by functions in ${\rm P}_g$. The trick is to decompose $g^n$ into its continuous $(g^n)^C$ and jump part $(g^n)^B$ and approximate each function independently.

 Let $n \in \mathbb N$. Observe that $g^n \in \operatorname{AC}_g([a,b])$. For $x\in D_g$, we have
 \[
 (g^n)'_g(x) = \frac{g^n(x^+)-g^n(x)}{g(x^+)-g(x)} = \sum_{k=0}^{n-1} g(x^+)^k g(x)^{n-1-k},
 \]
 and hence $|(g^n)'_g(x)| \leq n \|g\|_{\infty}^{n-1}$. Therefore, $(g^n)'_g \in L^2_{g^B}([a,b])$. 
 Enumerate $D_g = \{x_k\}_{k\in\Lambda}$ with $\Lambda \subset \mathbb N$. We aim to show that 
 $f \cdot 1_{D_g}$, where $f\in L^2_{g^B}([a,b])$,
 belongs to the closure of ${\rm P}_g$ in $L^2_g([a,b])$. If $\Lambda$ is finite, the claim is holds since $f \cdot 1_{D_g}$ is just a finite linear combination of the functions $1_{\{x\}}$ with $x\in D_g$. Suppose instead that $\Lambda = \mathbb N$, then
 \[
 \| f \cdot 1_{D_g} - \sum_{k=1}^{m} f(x_k) 1_{\{x_k\}}\|^2_{L^2_g} = \int_{[a,b)} (f \cdot 1_{D_g \setminus \{x_1,\dots,x_m\}})^2\operatorname{d}\mu_g=\sum_{k=m+1}^{\infty} f(x_k)^2\Delta g(x_k) \xrightarrow[m\to \infty]{} 0.
 \]
 In particular, $(g^n)'_g 1_{D_g}$ belongs to the closure of ${\rm P}_g$ in $L^2_g([a,b])$.

 By \cite[Lemma 6.1]{VillaTojo}, we have the decomposition $g^n = (g^n)^C + (g^n)^B$, where
 \[
 (g^n)^B(x) = \int_{[a,x)}(g^n)'_g\operatorname{d}\mu_{g^B} = \int_{[a,x)}(g^n)'_g 1_{D_g}\operatorname{d}\mu_{g} = \sum_{t \in [a,x)}g^n(t^+)-g^n(t).
 \]
 Clearly, $(g^n)^B \in \operatorname{AC}_g([a,b]) \subset \operatorname{UC}_g([a,b])$ and, since we can approximate its $g$--derivative by $g$--polynomials in $L^2_g([a,b])$, we can approximate $(g^n)^B$ by $g$--polynomials in $\operatorname{UC}_g([a,b])$, recall the proof of Theorem \ref{densidadequiv}.

 Also, since $(g^n)^C = g^n - (g^n)^B$, it follows that $(g^n)^C \in \operatorname{AC}_g([a,b]) \subset \operatorname{UC}_g([a,b])$. Moreover,
 \begin{align*}
 (g^n)^C(x) &= (g^n)^C(a) + \int_{[a,x)} \left[(g^n)^C\right]'_g \operatorname{d}\mu_g = (g^n)^C(a)+\int_{[a,x)}(g^n)'_g - (g^n)'_g 1_{D_g} \operatorname{d}\mu_g \\ &= (g^n)^C(a)+\int_{[a,x)}(g^n)'_g \operatorname{d}\mu_{g^C},
 \end{align*}
 for all $x \in [a,b]$. Thus, $(g^n)^C \in \operatorname{AC}_{g^C}([a,b]) \subset \operatorname{UC}_{g^C}([a,b])$. Next, we show that the functions $g^C_{a,k}$ belong to the closure of ${\rm P}_g$ in $\operatorname{UC}_{g}([a,b])$. Whence, by Proposition~\ref{casotontogcontinuo} and the fact that $\operatorname{UC}_{g^C}([a,b]) \subset \operatorname{UC}_g([a,b])$ it will follow that $\operatorname{UC}_{g^C}([a,b])$ is contained in the closure of ${\rm P}_g$ in $\operatorname{UC}_g([a,b])$ proving we can uniformly approximate $(g^n)^C$ by $g$--polynomials.

 We proceed by induction. For $k=1$, let $\varepsilon > 0$ and let $p$ be a $g$--polynomial such that \[ {\|p-(1-1_{D_g})\|_{L^2_g}<\varepsilon}.\]  Then,
 \begin{align*}
 \left|\int_{[a,x)}p\operatorname{d}\mu_g - g^C_{a,1}(x)\right| & = \left|\int_{[a,x)}p\operatorname{d}\mu_g - \int_{[a,x)}1 \operatorname{d}\mu_{g^C}\right| \\
 & = \left|\int_{[a,x)}p - (1-1_{D_g}) \operatorname{d}\mu_{g}\right| \leq \int_{[a,x)}\left|p - (1-1_{D_g})\right| \operatorname{d}\mu_{g} \\
 & \leq \int_{[a,b)}\left|p - (1-1_{D_g})\right| \operatorname{d}\mu_{g} \leq \|1\|_{L^2_g} \|p - (1-1_{D_g})\|_{L^2_g}<\|1\|_{L^2_g} \varepsilon
 \end{align*}
 for all $x \in [a,b]$. Note that the function $x\mapsto \int_{[a,x)}p\operatorname{d}\mu_g$ is itself a $g$--polynomial.

 Assume now the induction hypothesis. In particular, we have then that $g^C_{a,k}$ belongs to the closure of ${\rm P}_g$ in $L^2_{g}([a,b])$. Let $\varepsilon > 0$ and let $p$ be a $g$--polynomial such that $\|p-(k+1)g^C_{a,k}(1-1_{D_g})\|_{L^2}<\varepsilon$. Proceeding as before, we obtain
 \[
 \left|\int_{[a,x)}p\operatorname{d}\mu_g - g^C_{a,k+1}(x)\right| \leq \|1\|_{L^2} \varepsilon,
 \]
 for all $x\in[a,b]$. \qedhere
\end{proof}

\begin{rem}
 It turns out that the functions $p\circ g$ where $p$ is a polynomial are always dense in $\operatorname{UC}_g([a,b])$, independently of the continuity of $g$ as we proved on Proposition~\ref{casotontogcontinuo}. The idea is to use again Corollary~\ref{isometricgimagen} together with the Weierstrass approximation theorem. In this setting, however, $g([a,b])$ is not necessarily equal to $[g(a),g(b)]$, so an extension theorem such as the one proved in \cite[Theorem 3.13 and Remark 3.15]{VillaTojo} is needed.
\end{rem}

As a consequence of the previous theorem, we obtain the following corollary which reduces the Weierstrass approximation problem to the study of a family of functions.

\begin{cor}\label{corequiv}
 Let $g:[a,b]\to \mathbb R$ be a derivator. The following are equivalent:
 \begin{enumerate}
 \item[(1)] ${\rm P}_g$ is dense in $\operatorname{UC}_g([a,b])$.
 \item[(2)] For an arbitrary $x_0 \in D_g$, the bounded and nonincreasing sequence
 \[
 \(\frac{\operatorname{gram}(1,g_{x_0,1},\dots,g_{x_0,k})}{\operatorname{gram}(g_{x_0,1},\dots,g_{x_0,k})}\)_{k\in\bN}\]
converges to $\Delta g(x_0)$.
 \item[(3)] For an arbitrary $x_0 \in [a,b]$, the bounded and nonincreasing sequence
 \[
 \(\frac{\operatorname{gram}(1,g_{x_0,1},\dots,g_{x_0,k})}{\operatorname{gram}(g_{x_0,1},\dots,g_{x_0,k})}\)_{k\in\bN}
 \]
 converges to $ \Delta g(\beta)$ where $\beta=\max g^{-1}(\{g(x_0)\})$.
 \item[(4)] For an arbitrary $x_0\in [a,b]$, the function $1-1_{\{\beta\}}$ belongs to the closure of $\operatorname*{span}\{g_{x_0,k}\}_{k=1}^\infty$ in $L^2_g([a,b])$ where $\beta=\max g^{-1}(\{g(x_0)\})$.
 \end{enumerate}
\end{cor}
\begin{proof}
 By Theorem~\ref{equivalenciagpolinomiosdensidad}, ${\rm P}_g$ is dense in $\operatorname{UC}_g([a,b])$ if and only if the indicator functions $1_{\{x_0\}}$ belong to the closure of ${\rm P}_g$ in $L^2_g([a,b])$, for $x_0 \in D_g$. This holds if and only if, by Lemma~\ref{distanciaprodescalar},
 \[
 \frac{\operatorname{gram}(1_{\{x_0\}},1,g_{x_0,1},\dots,g_{x_0,k})}{\operatorname{gram}(1,g_{x_0,1},\dots,g_{x_0,k})}\xrightarrow[k\to\infty]{}0.
 \]
 Define ${\rm P}_k: = \operatorname{span}\{1,g_{x_0,1},\dots,g_{x_0,k}\}\subset \operatorname*{L}^2_g([a,b])$. Denote $p$ the orthogonal projection of the function $1_{\{x_0\}}$ over ${\rm P}_k$, so
 \[
 p = \sum_{m=0}^{k} a_mg_{x_0,m} \text{ and } 1_{\{x_0\}}-p\in {\rm P}_k^\perp.
 \]
 We have that
 \[
 d(1_{\{x_0\}},\operatorname{P}_k)^2 = \langle 1_{\{x_0\}},1_{\{x_0\}}\rangle_{L^2_g} - \langle p,1_{\{x_0\}}\rangle_{L^2_g} = \Delta g(x_0)-\Delta g(x_0)a_0.
 \]
 However,
 \[
 \operatorname*{Gram}(1,g_{x_0,1},\dots,g_{x_0,k})\begin{pmatrix}
	a_0\\a_1\\\vdots\\a_k
 \end{pmatrix} = \begin{pmatrix}
	\langle 1_{\{x_0\}},1\rangle_{L^2_g} \\\langle 1_{\{x_0\}},g_{x_0,1}\rangle_{L^2_g}\\\vdots\\\langle 1_{\{x_0\}},g_{x_0,k}\rangle_{L^2_g}
 \end{pmatrix} = \begin{pmatrix}
	\Delta g(x_0) \\0\\\vdots\\0
 \end{pmatrix}.
 \]
 Hence, by Cramer's rule,
 \[
 d(1_{\{x_0\}},\operatorname{P}_k)^2 = \Delta g(x_0) - \Delta g(x_0)^2 \operatorname{gram}(g_{x_0,1},\dots,g_{x_0,k}) \operatorname{gram}(1,g_{x_0,1},\dots,g_{x_0,k})^{-1}.
 \]
 The above tends to 0 if and only if
 \[
 \frac{\operatorname{gram}(1,g_{x_0,1},\dots,g_{x_0,k})}{\operatorname{gram}(g_{x_0,1},\dots,g_{x_0,k})} \xrightarrow[k\to\infty]{} \Delta g(x_0).
 \]
 Thus, (1) and (2) are equivalent. Take an arbitrary $x_0\in[a,b]$. Clearly, (3) implies (1) and (2). Let $\beta$ denote the maximum of the set $g^{-1}(\{g(x_0)\})$ and note that $g$ is constant on the interval $[x_0,\beta]$. Thanks to \cite[Proposition 3.14]{Cora2023}, the following formula holds for all $x\in\mathbb [a,b]$
 \[
 g_{x_0,n}(x) = \sum_{k=0}^{n} {n \choose k} g_{x_0,k}(\beta)g_{\beta,n-k}(x).
 \]
The functions $g_{x_0,k}$ are constant on $[x_0,\beta]$ so $g_{x_0,k}(\beta)=0$ for $k\geq1$. Thus, $g_{x_0,n}=g_{\beta,n}$ for all $n$. Therefore, if $\beta\in D_g$, (3) follows from (2). If $\beta\notin D_g$ then, by Lemma~\ref{lemmaintervalos}, $g^{-1}(\{g(x_0)\})$ has null $g$--measure.

 We show that if ${\rm P}_g$ is dense in $\operatorname{UC}_g([a,b])$ then the vector space $\operatorname{span}\{g_{x_0,k}\}_{k=1}^\infty$ is dense on the closed subspace
 \[
 \{f\in \operatorname{UC}_g([a,b]) \,: \,f(x_0)= 0\}
 \]
 on $\operatorname{UC}_g([a,b])$ for all $x_0\in[a,b]$. Indeed, fix $x_0\in[a,b]$ and let $f\in \operatorname{UC}_g([a,b])$ such that $f(x_0)=0$. Take any $\varepsilon>0$, there exists $p\in{\rm P}_g$ such that
 \[
 |p(x)-f(x)|<\frac{\varepsilon}{2}
 \]
 for all $x\in[a,b]$. The $g$--polynomial $p$ admits an expression 
 \[
 p(x) = \sum_{j\geq 0} \lambda_j g_{x_0,j}(x)
 \]
 for all $x\in[a,b]$ with $\lambda_j\in\mathbb F$ all zero except for a finite amount. Thus, 
 \[
 |\lambda_0 | = \left|\sum_{j\geq 0} \lambda_j g_{x_0,j}(x_0)-f(x_0)\right| = |p(x_0)-f(x_0)|<\frac{\varepsilon}{2}.
 \]
Hence,
\[
\left|\sum_{j\geq 1} \lambda_j g_{x_0,j}(x)-f(x)\right| = \left|\sum_{j\geq 0} \lambda_j g_{x_0,j}(x)-f(x)-\lambda_0\right|\leq |p(x)-f(x)|+|\lambda_0|<\varepsilon
\] 
for all $x\in[a,b]$.

 So if (1) holds and $\beta\notin D_g$, by Proposition~\ref{density}, $\{f\in \operatorname{UC}_g([a,b]) \,: \,f(x_0)= 0\}$ is dense on $L^2_g([a,b])$. Thus,
 $$1\in \overline{\operatorname{span}\{g_{x_0,k}\}_{k=1}^\infty} \text{ in } L^2_g([a,b])$$
and by Lemma~\ref{distanciaprodescalar} (3) follows. 

Assume (4) now. For all $x_0\in D_g$, $1-1_{\{x_0\}}\in\overline{\operatorname{span}\{g_{x_0,k}\}_{k=1}^\infty}$. Hence, $1_{\{x_0\}}\in\overline{\operatorname{span}\{g_{x_0,k}\}_{k=0}^\infty}=\overline{\rm P}_g$ and by Theorem~\ref{equivalenciagpolinomiosdensidad} (1) holds.

Take an arbitrary $x_0\in[a,b]$. Define $\beta$ as the maximum of the set $g^{-1}(\{g(x_0)\})$. If $\beta\notin D_g$, $1_{\{\beta\}}=0$ on $L^2_g([a,b])$ and (4) follows from (3). If $\beta\notin D_g$, consider $p\in \operatorname{span}\{g_{x_0,k}\}_{k=1}^\infty$. Note that the function $p-(1-1_{\{\beta\}})$ vanishes at $\beta$, as we showed above $g_{x_0,k}=g_{\beta,k}$, for all $k$, hence,
\begin{align*}
	||p-(1-1_{\{\beta\}})||^2_{L^2_g} & = \int_{[a,b)} (p-(1-1_{\{\beta\}}))^2\operatorname*{d}\mu_g
	= \int_{[a,b)\setminus\{\beta\}} (p-1)^2\operatorname*{d}\mu_g\\
	&= \int_{[a,b)} (p-1)^2\operatorname*{d}\mu_g-\Delta g(\beta)\\
	&=||p-1||^2_{L^2_g} -\Delta g(\beta).
\end{align*}
By (3) we can make the above as small as we want.
\end{proof}

Although the asymptotic behavior of this type of matrices has been studied—such as that of Gram matrices, or even Hankel matrices \cite{Min2021a,Bogatskiy2016,Its_2011} (the Gram matrix of classical monomials is a Hankel matrix)—the particular matrices that arise here in the context of Stieltjes calculus still present certain challenges.

Nevertheless, some results can still be stated about these determinants. Given  that 
\[
g_{a,n}(x) = \sum_{k=0}^{n} {n \choose k} g_{a,k}(x_0)g_{x_0,n-k}(x),
\]
for all $x\in[a,b]$---see \cite[Proposition 3.14]{Cora2023}, we can obtain the coefficients to change basis on $\operatorname{P}_k$, the space of $g$--polynomials with degree at most $k\in\mathbb N$, from that of  $g$--polynomials centered at $a$ to that of $g$--polynomials centered at $x_0$. Denote by $C$ the associated change of basis matrix. Hence,
\[
 \operatorname{Gram}(1,g_{a,1},\dots,g_{a,k}) = C^T \operatorname{Gram}(1,g_{x_0,1},\dots,g_{x_0,k}) C.
\]
Since $C$ is a triangular matrix with ones in its diagonal we have that
\[
 \operatorname{gram}(1,g_{a,1},\dots,g_{a,k}) = \operatorname{gram}(1,g_{x_0,1},\dots,g_{x_0,k}).
\]
Suppose $g$ is continuous and let $i,j\in\{1,\dots,k+1\}$. Then,
\begin{align*}
\operatorname{Gram}(1,g_{a,1},\dots,g_{a,k})_{ij} &= \int_{[a,b)} g_{a,i-1}g_{a,j-1}\operatorname*{d}\mu_g= \int_{[a,b)} g_{a,i+j-2}\operatorname*{d}\mu_g	\\
 &= \int_{[a,b)} g_{a,i+j-2}\operatorname*{d}\mu_g =\frac{g_{a,i+j-1}(b)-g_{a,i+j-1}(a)}{i+j-1}\\ 
 &= \frac{(g(b)-g(a))^{i+j-1}}{i+j-1}.
\end{align*}
Whence we obtain
\[
 \operatorname{gram}(1,g_{x_0,1},\dots,g_{x_0,k}) = (g(b)-g(a))^{(k+1)^2} H_{k+1}
\]
for all $x_0\in [a,b]$. Where $H_k$ denotes the determinant of the classical Hilbert matrix
\[
H_k = \begin{vmatrix}
 1 &\frac{1}{2} &\dots & \frac{1}{k}\\
 \frac{1}{2} &\frac{1}{3} & &\frac{1}{k+1}\\
 \vdots &\vdots & &\vdots \\
 \frac{1}{k} &\frac{1}{k+1} &\dots& \frac{1}{2k-1}\\
\end{vmatrix}.
\]

Although understanding the asymptotic behavior of these determinants is equivalent to proving the Weierstrass approximation theorem, the discontinuities in the derivator make their computation highly challenging. We proved Corollary~\ref{corequiv} because future progress via this approach cannot be ruled out. From now on, we will approach the problem from a different point of view.

\section{Auxiliary results}
\label{seccionaux}

In order to prove the Weierstrass Approximation Theorem for derivators with a finite amount of discontinuity points, a series of technical results is required. We have decided to write them separately so that they do not interfere with the exposition.

We start with the following definition, we define what it means for a function to be
$m$-times continuously differentiable when its domain is a degenerate interval. This ensures consistency for our results for all types of derivators.

\begin{dfn}
Let $\alpha\leq\beta$ be real numbers. We define
\[
 \operatorname{C}^\infty ([\alpha,\beta]):= \begin{dcases}
 \{f:[\alpha,\beta]\to \mathbb F\,|\, f\text{ is infinitely differentiable}\},& \text{if }\alpha < \beta,\\
 \mathbb F^\infty \equiv \prod_{n\in\mathbb N}\mathbb F,& \text{if }\alpha = \beta.
 \end{dcases}
\]
Thus, a function $f:\{\alpha\}\to \mathbb F$ belongs to $f\in\operatorname{C}^\infty (\{\alpha\})$ if it is associated with a sequence $(a_n)_{n\in\mathbb N}\subset \mathbb F$. We will use the notation
\[
f^{(k)}(\alpha)= a_k,\quad k=0,1,2,\dots
\]
For a given nonnegative integer $m$, define
\[
 \operatorname{C}^m ([\alpha,\beta]):= \begin{dcases}
 \{f:[\alpha,\beta]\to \mathbb F\,|\, f\text{ is }m\text{-times continuously differentiable}\},& \text{if }\alpha < \beta,\\
 \mathbb F^{m+1},& \text{if }\alpha = \beta.
 \end{dcases}
\]
Again, a function $f:\{\alpha\}\to \mathbb F$ belongs to $f\in\operatorname{C}^m (\{\alpha\})$ if it is associated with a vector $(a_j)_{j=0}^m\in\mathbb F^{m+1}$. We will use the notation
\[
f^{(k)}(\alpha) = a_k,\quad k=0,1,\dots,m.
\]
\end{dfn}

\begin{lem}
 \label{extensioninterpoladora}
 Let $[a,b]\subset R$ be a nondegenerate interval, $n\in \mathbb N$ and two vectors
\[
\mathbf v=(v_1,v_2,\dots,v_n),\,\mathbf w=(w_1,w_2,\dots,w_n)\in \mathbb R^n.
\]
 Denote $M$ and $m$ the maximum and minimum of $\{v_1,w_1\}$ respectively. Fix $\varepsilon>0$.
 There exists an infinitely differentiable
 function $F:[a,b]\to (m-\varepsilon,M+\varepsilon)$ such that
 \begin{equation}\label{probleminter}
 \begin{array}{l}
 F^{(k-1)}(a)=v_k,\\
 F^{(k-1)}(b)=w_k,
 \end{array}\quad \text{ for all }k=1,\dots,n.
 \end{equation}
\end{lem}
\begin{proof}
Let $p:[a,b]\to \mathbb R$ be the unique polynomial with degree at most $2n-1$ that solves the Hermite interpolation problem
\[
 \begin{array}{l}
 p^{(k-1)}(a)=v_k,\\ p^{(k-1)}(b)=w_k,
\end{array}\quad \text{ for all }k=1,\dots,n.
\]
By continuity of $p$, there exists $\delta>0$ such that \[ p(x)\in(m-\varepsilon,M+\varepsilon) \text{ for all } x\in[a,a+\delta]\cup [b-\delta,b].\]
Let $\varphi:[a,b]\to [0,1]$ be a smooth bump function that satisfies:
\begin{enumerate}[itemsep=0mm]
 \item $\varphi(x)\equiv 0$ on $[a,a+\frac{\delta}{2}]\cup [b-\frac{\delta}{2},b]$.
 \item $\varphi(x)\equiv 1$ on $[a+\delta,b-\delta]$.
\end{enumerate}
Define, for $x\in[a,b]$,
\[
F(x)=(1-\varphi(x))p(x)+\varphi(x)\frac{M+m}{2}.
\]
Note that $F$ is a convex combination of the polynomial $p$ and the mid--point of the interval $[m,M]$. Clearly $F\in\operatorname*{C}^\infty([a,b])$. Also,
\[
\begin{array}{l}
 F(x)= p(x) \text{ if } x\in[a,a+\frac{\delta}{2}]\cup [b-\frac{\delta}{2},b],\\
 F (x)= \frac{M+m}{2} \text{ if } x\in[a+\delta,b-\delta].
\end{array}
\]
This means $F$ satisfies equation~\eqref{probleminter}. Besides, for $x\in[a+\frac{\delta}{2},a+\delta]\cup [b-\delta,b-\frac{\delta}{2}]$, $F(x)$ must belong to $(m-\varepsilon,M+\varepsilon)$
since this set is convex and both $p(x)$ and $(M+m)/2$ belong to it. Hence, $F$ satisfies the conditions of the theorem.\qedhere
\end{proof}

In the proof of Lemma~\ref{extensioninterpoladora} we look for a smooth function of which we have control on a finite number of its derivatives on the end points of the interval and which we can fit on the rectangle
$[a,b]\times(m-\varepsilon,M+\varepsilon)$ for any given $\varepsilon>0$.
The following Corollary is a straightforward consequence of Lemma~\ref{extensioninterpoladora}.

\begin{cor}
 \label{interpolarcota}
 Let $\alpha\leq\beta <\gamma\leq\delta$ be real numbers and $m$ nonnegative integer. Consider $f_1\in\operatorname{C}^m([\alpha,\beta])$ and $f_2\in\operatorname{C}^m([\gamma,\delta])$.
 For any $M>\max\{\|f_1\|_\infty,\|f_2\|_{\infty}\}$, there exists $F\in\operatorname{C}^m([\alpha,\delta])$ such that $f_1^{(k)}(\beta) = F^{(k)}(\beta)$ and $f_2^{(k)}(\gamma)=F^{(k)}(\gamma)$ for $k=0,1,\dots,m$
 and
 \[
 \begin{array}{l}
 F(x)= f_1(x), \text{ if } x\in[\alpha,\beta],\\
 F (x)= f_2(x), \text{ if } x\in[\gamma,\delta],
 \end{array}
 \]
 with $\|F\|_\infty < M$.
\end{cor}

\begin{lem}
 \label{pegado}
Let $n\in\mathbb N$. Consider $x_0\leq x_1\leq \dots\leq x_n$, $n+1$ points in $\mathbb R$. For $j\in\{1,\dots,n\}$, let $\widehat f_j:[x_{j-1},x_j]\to \mathbb F$ be such that $\widehat f_j\in\operatorname{C}^\infty ([x_{j-1},x_j])$. Assume the following holds,
\[
\widehat f_i^{(k)}(x_{i})=\widehat f_{i+1}^{(k)}(x_{i}) \text{ for all } k\leq i-1 \text{ and all }i\in\{1,\dots,n-1\}.
\]
%
For any fixed $\delta>0$, there exists a function $F\in \operatorname{C}^{n-1}([x_0,x_n])$ such that, for all $k=0,\dots,n-1$,
 \[
 |F^{(k)}(x)-\widehat f_j^{(k)}(x)|<\delta\text{ for }x\in[x_{j-1},x_j] \text{ and }j\in\{k+1,\dots,n\}.
 \]
\end{lem}
\begin{proof}

We prove this by induction on $n$. For $n=1$, the function $F=\widehat f_1$ verifies the assertion. Fix $n\geq 2$ and $\delta >0$. We have two cases.

\textbullet If $x_0\neq x_1$, let $\varepsilon>0$ be such that
\[
\varepsilon<\begin{dcases}
 \operatorname*{min}\left\{\frac{\delta}{2(x_n-x_1)},\delta\right\}, &\text{if } x_n\neq x_1,\\
 \delta, & \text{if } x_n=x_1.
\end{dcases}
\]

By the induction hypotheses, there exists some function $H\in\operatorname{C}^{n-2}([x_1,x_n])$ such that, for all $k=0,\dots,n-2$,
\[
 |H^{(k)}(x)-\widehat f_j^{(k+1)}(x)|<\varepsilon\text{ for }x\in[x_{j-1},x_j] \text{ and }j\in\{k+2,\dots,n\}.
\]
Fix a positive constant $M>\max\{\|H\|_{\infty},\|\widehat f_1'\|_{\infty}\}$. For $k\in\mathbb N$, $k>1/(x_1-x_0)$, by Corollary~\ref{interpolarcota}, there exists a function $H_k\in\operatorname*{C}^{n-2}([x_0,x_n])$ such that
\[
H_k(x)=\begin{cases}
 \widehat f_1'(x) & \text{ if }x\in \left[x_0,x_1-\frac{1}{k}\right],\\
 H(x) & \text{ if } x\in[x_1,x_n],
\end{cases}
\]
with $|H_k(x)|<M$ for all $x\in[x_0,x_n]$. Denote
\[
 F_k(x) = \widehat f_1(x_0) + \int_{x_0}^{x} H_k(t)\operatorname*{d}t, \quad x\in[x_0,x_n].
\]
The following holds,
\begin{enumerate}[itemsep=0mm]
 \item[(1)] $F_k\in\operatorname{C}^{n-1}([x_0,x_n])$,
 \item[(2)] For $k=1,\dots,n-1$, since $F_k'=H$ on $[x_1,x_n]$,
 \[
 |F^{(k)}(x)-\widehat f_j^{(k)}(x)|<\delta\text{ for }x\in[x_{j-1},x_j] \text{ and }j\in\{k+1,\dots,n\}.
 \]
 \item[(3)] For $j\in\{1,\dots,n\}$ and $x\in[x_{j-1},x_j]$,
 \begin{align*}
 |F_k(x)-\widehat f_j(x)|&= \left| F_k(x_{j-1}) + \int_{x_{j-1}}^{x} H_k(t)\operatorname*{d}t - \widehat f_j (x_{j-1}) - \int_{x_{j-1}}^{x}\widehat f_j'(t)\operatorname*{d}t \right|\\
 &\leq |F_k(x_{j-1})-\widehat f_j(x_{j-1})|+ \int_{x_{j-1}}^{x}|H_k(t)-\widehat f_j'(t)|\operatorname*{d}t\\
 &\leq |F_k(x_{j-2})-\widehat f_j(x_{j-2})|+ \int_{x_{j-2}}^{x_{j-1}}|H_k(t)-\widehat f_{j-1}'(t)|\operatorname*{d}t+\int_{x_{j-1}}^{x}|H_k(t)-\widehat f_j'(t)|\operatorname*{d}t\\
 &\leq \sum_{l=1}^{j} \int_{x_{l-1}}^{x_{l}}|H_k(t)-\widehat f_{j-1}'(t)|\operatorname*{d}t\leq \sum_{l=1}^{n} \int_{x_{l-1}}^{x_{l}}|H_k(t)-\widehat f_{j-1}'(t)|\operatorname*{d}t\\
 &\leq \int_{x_{1}-\frac{1}{k}}^{x_{1}}|H_k(t)-\widehat f_{1}'(t)|\operatorname*{d}t+\sum_{l=2}^{n} \int_{x_{l-1}}^{x_{l}}|H(t)-\widehat f_{j-1}'(t)|\operatorname*{d}t\\
 &\leq \frac{2M}{k}+\frac{\delta}{2}.
 \end{align*}
\end{enumerate}
Take any $k>\frac{4M}{\delta}$. Then $F=F_k$ is the function we looked for.
%

\textbullet If $x_0=x_1$ let $\varepsilon>0$ be such that
\[
\varepsilon<\begin{dcases}
 \operatorname*{min}\left\{\frac{\delta}{x_n-x_1},\delta\right\} &\text{if } x_n\neq x_1, \\
 \delta & \text{if } x_n=x_1.
\end{dcases}
\]

Again, $\widehat F'$ satisfies the induction hypotheses on $[x_1,x_n]$. Let $H\in\mathcal{C}^{n-2}([x_1,x_n])$ such that, for $k=0,\dots,n-2$,
\[
 |H^{(k)}(x)-\widehat f_j^{(k+1)}(x)|<\varepsilon\text{ for }x\in[x_{j-1},x_j] \text{ and }j\in\{k+2,\dots,n\}.
\]
Define the function,
\[
 F(x) = \widehat f_1(x_0) + \int_{x_0}^{x} H(t)\operatorname*{d}t, \quad x\in[x_0,x_n].
\]
Then, $F'(x)=H(x)$ for all $x\in[x_0,x_n]$, $F\in\operatorname*{C}^{n-1}([x_0,x_n])$ and $F(x_0)=\widehat f_1(x_0)=\widehat f_2(x_0)$. Also, for $j\in\{1,\dots,n\}$ and $x\in[x_{j-1},x_j]$,
\begin{align*}
 |F(x)-\widehat f_j(x)|&\leq \sum_{l=2}^{n} \int_{x_{l-1}}^{x_{l}}|H(t)-\widehat f_{j-1}'(t)|\operatorname*{d}t\leq \delta.
\end{align*}
$F$ is the function looked for.
\end{proof}

\begin{rem}
 Applying Lemma~\ref{pegado} recursively for a decreasing sequence of deltas that converge to $0$, for example, $\delta_m=1/m$, we get a sequence of functions $(F_m)_{m\in\mathbb N}\subset\operatorname{C}^{n-1}([x_0,x_n])$ such that
 $(F_m)_{m\in\mathbb N}$ converges uniformly on $[x_0,x_n]$ to
 \[
 G^0(x)=\widehat f_j(x),\text{ if } x\in[x_{j-1},x_j] \text{ with } j=1,\dots,n;
 \]
 $(F_m')_{m\in\mathbb N}$ converges uniformly on $[x_1,x_n]$ to
 \[
 G^1(x)=\widehat f_j'(x),\text{ if } x\in[x_{j-1},x_j] \text{ with } j=2,\dots,n;
 \]
 and, in general, $(F_m^{(k)})_{m\in\mathbb N}$ converges uniformly on $[x_k,x_n]$ to
 \[
 G^k(x)=\widehat f_j^{(k)}(x),\text{ if } x\in[x_{j-1},x_j] \text{ with } j=k+1,\dots,n,
 \]
 for $k=0,\dots,n-1$. Note that $G^k$ is continuous on $[x_k,x_n]$.
\end{rem}

\begin{lem}
 \label{weierdife}
 Fix $\varepsilon>0$.
 Let $n\in\mathbb N$. Consider $x_0\leq x_1\leq \dots\leq x_n$, $n+1$ points in $\mathbb R$. Let $ f_j\in\operatorname{C}^\infty ([x_{j-1},x_j])$ with $j\in\{1,\dots,n\}$. For any constants $c_k^j\in\mathbb R$ with $j=0,\dots,n-1$, $k=0,\dots,j$ and $c^j_0=1$, there exists a polynomial $p$ such that
\[
 \left|\sum_{k=0}^{j-1} c^{j-1}_k p^{(k)}(x) - f_j(x)\right|<\varepsilon \text{ for all }x\in[x_{j-1},x_j],\quad j=1,\dots,n.
\]
\end{lem}
\begin{proof}
Fix $\varepsilon>0$. Let $\widehat f_1 = f_1$. Take now $\widehat f_2$ as the unique solution of the differential equation with initial conditions
\[
\widehat f_2+c_1^1 \widehat f_2'=f_2 \text{ on } [x_1,x_2],\quad \widehat f_2(x_1)=\widehat f_1(x_1).
\]
Inductively, for $j=2,\dots,n$, let $\widehat f_j$ be the unique solution of
\[
\widehat f_j + c_1^{j-1}\widehat f_j'+\dots +c_{j-1}^{j-1}\widehat f_j^{(j-1)} = f_j\text{ on }[x_{j-1},x_j],\quad \widehat f_j^{(k)}(x_{j-1})=\widehat f_{j-1}^{(k)}(x_{j-1}) \text{ for all } k\leq j-2.
\]
By construction, the functions $\widehat f_j$ satisfy the hypotheses of Lemma~\ref{pegado}. Let $M$ be a positive constant greater than $|c_k^j|$, $j\in\{0,\dots,n-1\}$, $k\in\{0,\dots,j\}$.
Now, by Lemma~\ref{pegado}, there exists a function $F\in\operatorname{C}^{n-1}([x_0,x_n])$ such that, for all $k=0,\dots,n-1$,
\[
|F^{(k)}(x)-\widehat f_j^{(k)}(x)|<\frac{\varepsilon}{2nM}\text{ for }x\in[x_{j-1},x_j] \text{ and }j\in\{k+1,\dots,n\}.
\]
By the classical Weierstrass Approximation Theorem, let $p$ be a polynomial such that
\[
 |p^{(k)}(x)-F^{(k)}(x)|<\frac{\varepsilon}{2nM} \text{ for }x\in[x_0,x_n]\text{ and } k=0,\dots,n-1.
\]
Now, for $j=1,\dots,n$ and $x\in[x_{j-1},x_j]$,
\begin{align*}
 \left|\sum_{k=0}^{j-1} c^{j-1}_k p^{(k)}(x) - f_j(x)\right|&=\left|\sum_{k=0}^{j-1} c^{j-1}_k p^{(k)}(x) - \sum_{k=0}^{j-1} c^{j-1}_k \widehat f_j^{(k)}(x)\right|\\
 &\leq\sum_{k=0}^{j-1} |c^{j-1}_k\| p^{(k)}(x)-\widehat f_j^{(k)}(x)|\\
 &\leq\sum_{k=0}^{j-1} |c^{j-1}_k\| p^{(k)}(x)-F^{(k)}(x)|+\sum_{k=0}^{j-1} |c^{j-1}_k\|F^{(k)}(x)-\widehat f_j^{(k)}(x)|< \varepsilon.
\end{align*}
\end{proof}

From Lemma~\ref{weierdife} we derive the following corollary.
\begin{cor}
 Fix $n\in\mathbb N$ and consider
 $x_0\leq x_1\leq \dots\leq x_n$,
 $n+1$ points in $\mathbb R$. Let $ f_j:[x_{j-1},x_j]\to \mathbb F$, $ f_j\in\operatorname{C}^\infty ([x_{j-1},x_j])$ with $j\in\{1,\dots,n\}$. For any positive constants $c_k^j\in\mathbb R$ with $j=0,\dots,n-1$, $k=0,\dots,j$ and $c^j_0=1$. There exists a sequence of polynomials $p_m$ such that
 \begin{align*}
 &p_m \rightarrow f_1\quad \text{uniformly on }[x_0,x_1],\\
 &p_m + c_1^1 p_m' \rightarrow f_2\quad \text{uniformly on }[x_1,x_2],\\
 &\quad\quad\quad\quad\quad\quad\vdots\\
 &p_m + c_1^{n-1} p_m'+ \dots + c_{n-1}^{n-1}p_m^{(n-1)} \rightarrow f_n\quad \text{uniformly on }[x_{n-1},x_n].
 \end{align*}
\end{cor}

\section{Weierstrass Approximation Theorem for Stieltjes Calculus}
\label{pruebaseccion}

We now have the necessary results to prove the main theorem of this work. This entire section is devoted to prove it. We state it below.

\begin{thm}[Weierstrass Approximation Theorem for Stieltjes Calculus]
 \label{WAT}
 Let $[a,b]$ be an interval and $g$ a derivator with a finite amount of discontinuity points on $[a,b]$. Then, ${\rm P}_g$ is dense in $\operatorname*{UC}_g([a,b])$.
\end{thm}

In particular, following Theorem~\ref{WAT} and \cite[Section 3]{Cora2023}, infinitely $g$--differentiable functions are dense on the space of uniformly $g$--continuous functions.

From now on, $g$ will denote a derivator defined on $[a,b]$ with a finite amount of discontinuity points.
Let $\{x_1,\dots,x_{n-1}\}=D_g\cap[a,b)$ with $x_j < x_k$ for $j<k$ and denote $x_0=a$ and $x_n=b$.

\begin{rem}
 Note that $g$ and $g_{a,1}=g-g(a)$ are equivalent derivators. They define the same pseudometric, Lebesgue-Stieltjes integral and Stieltjes differentiation. Sets like $C_g$ and $D_g$ are also
 the same for both derivators. This means they define the same Stieltjes continuity (uniform or not) and Stieltjes monomials. Thus, proving Theorem~\ref{WAT} for $g_{a,1}$ is proving it for $g$ also.

 Hence, without loss of generality, we will assume all derivators (including $g^C$ and $g^B$) satisfy $g=g_{a,1}$, or what is equivalent, that $g(a)=0$.
\end{rem}

\begin{dfn}
Let $V_j$ be a normed linear space for $j=1,\dots,n$. We consider on $\bigoplus_{j=1}^n V_j$ the norm
\[
\|v_1+\cdots+v_n\|=\max\{\|v_j\|_{V_j}:j=1,\dots,n\},
\]
for $v_j\in V_j$.
\end{dfn}

\begin{lem}
 \label{psi}
The linear operator
\[
\begin{tikzcd}[row sep= 0em]
		\Psi:\operatorname{UC}_g([a,b]) \arrow{r} & \operatorname{UC}_{g^C}([x_{0},x_1])\oplus\bigoplus\limits_{j=2}^n \operatorname{UC}_{g^C}((x_{j-1},x_j]) \\
		f \arrow[mapsto]{r} & (\restr{f}{|[x_0,x_1]},\restr{f}{|(x_1,x_2]}\dots,\restr{f}{|(x_{n-1},x_n]}),
\end{tikzcd}
\]
is an isometric isomorphism. Where $\restr{f}{|(x_{j-1},x_j]}$ denotes the usual restriction.
\end{lem}
\begin{proof}
$\Psi$ is well-defined. Let $j\in\{2,\dots,n\}$. Take $f\in \operatorname{UC}_g([a,b])$, clearly, $f$ is uniformly $g$--continuous on $(x_{j-1},x_j]$. For any $x,y\in(x_{j-1},x_j]$, $x<y$,
\begin{align*}
g(y)-g(x)&=g^C(y)-g^C(x)+g^B(y)-g^B(x)\\[ 0.1em]
&=g^C(y)-g^C(x)+\sum\limits_{t\in D_g\cap[x,y)}\Delta g(t)\\[ -.6em]
&=g^C(y)-g^C(x).
\end{align*}
Therefore, $f$ is uniformly $g^C$--continuous on $(x_{j-1},x_j]$. For $j=1$, we have two cases: If $a\in D_g$, $x_1=x_0=a$ and $f_{|[x_0,x_1]}=f_{|\{a\}}\in \operatorname{UC}_{g^C}(\{a\})$.
If $a\notin D_g$, $g(y)-g(x)=g^C(y)-g^C(x)$ for all $x,y\in[a,x_1]$ with $x<y$ as before and $f_{|[a,x_1]}\in\operatorname{UC}_{g^C}([a,x_1])$.


Observe
\[ \sup_{x\in[a,b]}|f(x)|=\max\left\{\sup_{x\in (x_{j-1},x_j]} |f(x)|: 2\leq j \leq n\right\}\cup \left\{\sup_{x\in [x_{0},x_1]} |f(x)|\right\},\]
so, by the definition of the norms involved, $\Psi$ is an isometry. We now show that $\Psi$ is onto. Let $f_j\in\operatorname{UC}_{g^C}((x_{j-1},x_j])$ for $j=2,\dots,n$ and let $f_1\in\operatorname{UC}_{g^C}([x_0,x_1])$. Define
\[
f(x)=\begin{dcases}
f_1(x), & x\in[x_0,x_1],\\
f_j(x), & x\in(x_{j-1},x_j], \quad j=2,\dots,n.
\end{dcases}
\]
First, $f$ is uniformly $g$--continuous on $[a,b]$. Indeed, fix $\varepsilon>0$. There is $\delta>0$ such that
\begin{align*}
 & \text{ if } \,x,y\in [x_{0},x_1] \text{ with } |g^C(x)-g^C(y)|<\delta \Rightarrow |f_1(x)-f_1(y)|<\varepsilon,\\
 & \text{ if } \,x,y\in (x_{j-1},x_j] \text{ with } |g^C(x)-g^C(y)|<\delta \Rightarrow |f_j(x)-f_j(y)|<\varepsilon,\quad j\in\{2,\dots,n\}.
\end{align*}
Assume that $\delta <\min_{x\in D_g}\Delta g(x)$.
Take $x,y\in[a,b]$, $x<y$ and assume there is $j\in\{1,\dots,n\}$ such that $x\leq x_j < y$. Then,
\begin{align*}
 g(y)-g(x)&=g(y)-g(x_j^+)
 +g(x_j^+)-g(x_j)
 +g(x_j)-g(x)>\delta.
\end{align*}
Thus, if $g(y)-g(x)<\delta$, both points must be contained in one of the following disjoint intervals
\begin{align*}
 &[x_0,x_1],\\
 &(x_{j-1},x_{j}], \quad j\in\{2,\dots,n\},
\end{align*}
but on those intervals $f$ coincides with $f_j$ and $g$ with $g^C$ so $f$ is uniformly $g$--continuous and \[ \Psi(f)=(f_1,\dots,f_n).\qedhere\]
\end{proof}

The previous and next result show that a uniformly $g$--continuous function $f\in\operatorname{UC}_{g}([a,b])$ is a piecewise uniformly $g^C$--continuous function (and thus piecewise continuous) that may have
jumps on the points of discontinuity of $g$.

\begin{rem}
 \label{extensa}
 For $i=2,\dots,n$, the spaces $\operatorname{UC}_{g^C}((x_{i-1},x_i])$ and $\operatorname{UC}_{g^C}([x_{i-1},x_i])$ are isometrically isomorphic. Indeed, it is enough to map every $f\in\operatorname{UC}_{g^C}((x_{i-1},x_i])$ to $\widetilde f\in\operatorname{UC}_{g^C}([x_{i-1},x_i])$ where
\begin{equation*}
 \widetilde f(x)=\begin{dcases}
 f(x), & x\in (x_{i-1},x_i],\\
 \lim_{y\to x_{i-1}^+}f(y), & x=x_{i-1}.
 \end{dcases}
\end{equation*}
\end{rem}

\begin{lem}
 \label{philem}
 For $i=1,\dots,n$, the spaces $\operatorname{UC}_{g^C}([x_{i-1},x_i])$ and $\operatorname{C}([g^C(x_{i-1}),g^C(x_i)])$ are isometrically isomorphic.
\end{lem}
\begin{proof}
This result is a consequence of Corollary~\ref{isometricgimagen} together with the fact that, since $g^C$ is continuous and nondecreasing, $g^C([x_{i-1},x_i])=[g^C(x_{i-1}),g^C(x_i)]$. The rest follows from the fact that
continuity and uniform continuity are equivalent on compact sets.
\end{proof}

Recall Definition~\ref{gmonomio}.
Since for all $x_0\in [a,b]$, the set $\{g_{x_0,n}\}_{n=0}^\infty$ is a generating set of $\mathrm{P}_g$. From now on we consider all $g$--monomials centered at $x_0=a$ and, thus, we will use the notation $g_n\equiv g_{a,n}$, $g_n^C\equiv g^C_{a,n}$ and $g_n\equiv g^B_{a,n}$.

\begin{dfn}
We define
\[
\begin{tikzcd}[row sep= 0em]
 \Phi:\bigoplus\limits_{i=1}^n \operatorname{C}([g^C(x_{i-1}),g^C(x_i)]) \arrow{r}& \operatorname{UC}_{g^C}([x_{0},x_1])\oplus\bigoplus\limits_{i=2}^n \operatorname{UC}_{g^C}((x_{i-1},x_i]) \\
 (f_1,\dots,f_n) \arrow[mapsto]{r} & (f_1\cdot g^C, \dots,f_n\cdot g^C).
\end{tikzcd}
\]
\end{dfn}

Lemma~\ref{philem} together with Remark~\ref{extensa} show that $\Phi$ is an isometric isomorphism. This together with Lemma~\ref{psi} constructs the isometric isomorphism \[ \Phi^{-1}\Psi:\operatorname{UC}_g([a,b])\to\bigoplus_{i=1}^n \operatorname*{C}([g^C(x_{i-1}),g^C(x_i)]).\]

Hence, ${\rm P}_g$ is dense on $\operatorname*{UC}_g([a,b])$ if and only if  $\Phi^{-1}\Psi(\operatorname*{P}_g)$ is dense on $\bigoplus_{i=1}^n \operatorname*{C}([g^C(x_{i-1}),g^C(x_i)])$.
We will show exactly that, for which it is crucial to determine the structure of $\Phi^{-1}\Psi(\operatorname*{P}_g)$.

\begin{lem}
	\label{pgfactor}
 Let $p_g$ be an arbitrary $g$--polynomial. That is,
 \[
 p_g(x)=\sum_{k=0}^m \alpha_k g_k(x),
 \]
 for $x\in[a,b]$ and $\alpha_k\in \mathbb F$.
 Denote,
 \[
 p(x)=\sum_{k=0}^m \alpha_k x^k,
 \]
 for $x\in[a,b]$. Then,
 \begin{equation}
 \label{formulaphipsi}
 \Phi^{-1}\Psi(p_g)=\left(p,p+p'g_1^B(x_2),\dots,\sum_{k=0}^{n-1} p^{(k)}\frac{g_k^B(x_n)}{k!}\right).
 \end{equation}
\end{lem}
\begin{proof}
By Theorem~\ref{formulasmonomios} and Proposition~\ref{propiedadesmonomios}, for $x\in[a,b]$,
\begin{align*}
 p_g(x)&=\sum_{k=0}^m \alpha_k \sum_{i=0}^{k}\binom{k}{i}g^C_{k-i}(x)g^B_{i}(x)=\sum_{k=0}^m \alpha_k \sum_{i=0}^{k}\binom{k}{i}g^C(x)^{k-i}g^B_{i}(x)\\
 &=\sum_{k=0}^{m} p^{(k)}(g^C(x))\frac{g_k^B(x)}{k!}.
\end{align*}
Although in this case it follows from direct calculation, the formula
\[ p_g(x)=\sum_{k=0}^{m} p^{(k)}(g^C(x))\frac{g_k^B(x)}{k!}\]
is proven more generally in \cite[Proposition 4.3]{Cora2023}.

Besides, $g_k^B$ is uniformly $g^B$--continuous and hence it must be constant where $g^B$ is constant. Thus, since $g^B$ is constant on the intervals $[x_0,x_1]$ and $(x_{i-1},x_{i}]$ for $i=2,\dots,n$, so is $g^k_B$ for all $k$.
Furthermore, following Proposition~\ref{propiedadesmonomios}, the $g^B_k$ are all nonnegative on $[a,b]$ and
\[
 g^B_k(x)=0\text{ for }x\in[a,x_k],\ k=1,\dots,n.
\]
 We also have that $g^B_k(x)=0$ on $[a,b]$ if $k\geq n$. Therefore,
\[
\begin{array}{ll}
 \text{on } [x_0, x_1] & \text{all } g_k^B \text{ vanish}, \\
 \text{on } (x_1, x_2] & \text{all } g_k^B \text{ vanish except for } k = 1, \\
 \text{on } (x_2, x_3] & \text{all } g_k^B \text{ vanish except for } k = 1,2, \\
 & \hspace{1.5em} \vdots \\
 \text{on } (x_{n-1}, x_n] & \text{all } g_k^B \text{ vanish except for } k = 1, \dots, n-1.
\end{array}
\]
With the above, for $x\in[x_0,x_1]$,
\[
p_{g}(x)=p(g^C(x)).
\]
For $i\in\{2,\dots,n\}$ and $x\in(x_{i-1},x_{i}]$,
\[
p_g(x)=\sum_{k=0}^{i-1} p^{(k)}(g^C(x))\frac{g_k^B(x_i)}{k!}.
\]
Hence, $\Phi^{-1}\Psi$ maps $p_g$ to
\[
 \Phi^{-1}\Psi(p_g)=\left(p,p+p'g_1^B(x_2),\dots,\sum_{k=0}^{n-1} p^{(k)}\frac{g_k^B(x_n)}{k!}\right).
\qedhere\]
\end{proof}

We then have that, following equation~\eqref{formulaphipsi}, $\Phi^{-1}\Psi(\operatorname*{P}_g)$ is dense on $\bigoplus_{i=1}^n \operatorname*{C}([g^C(x_{i-1}),g^C(x_i)])$ if and only if for any $n$-tuple $(f_1,\dots,f_n)$ of continuous functions
there exist a sequence of polynomials $p_m$ such that
\begin{align}
 &p_m \rightarrow f_1\quad \text{uniformly on }[g^C(x_0),g^C(x_1)],\notag\\
 &p_m + g_1^B(x_2) p_m' \rightarrow f_2\quad \text{uniformly on }[g^C(x_1),g^C(x_2)],\notag\\
 & \hspace{8em} \vdots \notag\\[ 0em]
 &p_m + g_1^B(x_n) p_m' +\dots + \frac{g_{n-1}^B(x_n)}{(n-1)!} p_m^{(n-1)} \rightarrow f_n\quad \text{uniformly on }[g^C(x_{n-1}),g^C(x_n)].\label{eq:seqpol}
\end{align}
\begin{rem}
 First, note that $g^C(x_0)=g^C(a)=0$. Therefore, $g^C([a,b])=[0,g^C(b)]$. It could happen that $g^C([a,b])=\{0\}$. Whence,
 \[
 \operatorname{UC}_g([a,b])=\bigoplus_{i=1}^n \operatorname{C}(\{0\})=\mathbb F^n.
 \]
 In this case, Theorem~\ref{WAT} trivially holds since, for $p(x)=\sum_{k=0}^{n-1} a_k x^k$, clearly $p^{(k)}(0)=k! a_k$. Thus, problem~\eqref{eq:seqpol} transforms into
 \[
 \begin{bmatrix}
 1&0&0&\dots&0\\
 1&g_1^B(x_2)&0&\dots&0\\
 &&\vdots&&\\
 1&g_1^B(x_n)&g_2^B(x_n)&\dots&g_{n-1}^B(x_n)\\
 \end{bmatrix}
 \begin{bmatrix}
 a_0\\
 a_1\\
 \vdots\\
 a_{n-1}
 \end{bmatrix}
 =
 \begin{bmatrix}
 f_1(0)\\
 f_2(0)\\
 \vdots\\
 f_n(0)
 \end{bmatrix}
 \]
 which is clearly solvable. Notice it can happen that $g^C(x_{i-1})=g^C(x_{i})$ arbitrarily.
\end{rem}

\begin{thm}
 \label{weierstrass}
 Let $[0,b]$ be a given interval. Consider $n+1$ points $0=x_0\leq x_1\leq x_2\leq\dots\leq x_{n-1}\leq x_n= b$. For any $n$-tuple
 $(f_1,\dots,f_n)$ of continuous functions and any positive constants $c_k^i\in\mathbb R$ there exist a sequence of polynomials $p_m$ such that
\begin{align*}
 &p_m \rightarrow f_1\quad \text{uniformly on }[x_0,x_1],\\
 &p_m + c_1^1 p_m' \rightarrow f_2\quad \text{uniformly on }[x_1,x_2],\\
 &\quad\quad\quad\quad\quad\quad\vdots\\
 &p_m + c_1^{n-1} p_m'+ \dots + c_{n-1}^{n-1}p_m^{(n-1)} \rightarrow f_n\quad \text{uniformly on }[x_{n-1},x_n].
\end{align*}
\end{thm}
\begin{proof}
 Let $(f_1,\dots,f_n)$ be any $n$-tuple of continuous functions. We show that, given $\varepsilon>0$ there exist a polynomial $p$
 such that
 \[
 \left|\sum_{k=0}^{j-1} c^{j-1}_k p^{(k)}(x) - f_j(x)\right|<\varepsilon \text{ for all }x\in[x_{j-1},x_j]
 \]
 with $j\in\{1,\dots,n\}$. This clearly is equivalent to the thesis of the theorem. For $j=1,\dots,n$, let $h_j$ be infinitely differentiable functions such that
 \[
 |h_j(x)-f_j(x)|<\frac{\varepsilon}{2}\text{ for }x\in[x_{j-1},x_j].
 \]
 The functions $h_j$ are guaranteed to exist using Weierstrass' Approximation Theorem. Applying Lemma~\ref{weierdife} to the $n$-tuple $(h_1,\dots,h_n)$, there exist
 a polynomial $p$ such that
 \[
 \left|\sum_{k=0}^{j-1} c^{j-1}_k p^{(k)}(x) - f_j(x)\right|<\frac{\varepsilon}{2} \text{ for all }x\in[x_{j-1},x_j],
 \]
 with $j=1,\dots,n$. Therefore, for $j=1,\dots,n$ and $x\in[x_{j-1},x_j]$,
\[
 \left|\sum_{k=0}^{j-1} c^{j-1}_k p^{(k)}(x) - f_j(x)\right|\leq \left|\sum_{k=0}^{j-1} c^{j-1}_k p^{(k)}(x) - h_j(x)\right|+|h_j(x)-f_j(x)|\leqslant \frac{\varepsilon}{2}+\frac{\varepsilon}{2}=\varepsilon.\qedhere
\]
\end{proof}
Finally, Theorem~\ref{WAT} follows from Theorem~\ref{weierstrass}.
\begin{proof}[Proof of Theorem~\ref{WAT}]
 $\Phi^{-1}\Psi(\operatorname{P}_g)$ is a dense subset of $\bigoplus_{j=1}^n \operatorname{C}([g^C(x_{j-1}),g^C(x_j)])$, by Theorem~\ref{weierstrass}. Since
 $\Phi^{-1}\Psi$ is an isometric isomorphism, $\mathrm{P}_g$ must be a dense subset of $\operatorname{UC}_g([a,b])$.
\end{proof}

\section{Weierstrass for discontinuous derivators}
\label{secfinal}
We only proved the Weierstrass approximation theorem on derivators on which $D_g$ is finite. However, it is easy to extend the result to some class of derivators with infinite amount of discontinuities. Let $g:[a,b]\to \mathbb R$ be such that $g^C=0$, $D_g=\{x_n\}_{n\in\mathbb N}$ with $x_1=a$ and, for every $n\in\mathbb N$, $x_n < x_{n+1}$, that is $D_g$ is well-ordered. Then, $D_g$ looks like the

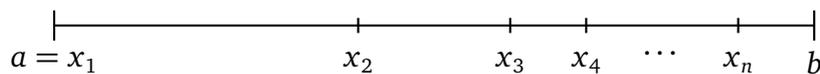
\begin{center}
 \begin{tikzpicture}
 \draw[thick] (0,0) -- (10,0);

 \draw[thick] (0,-0.2) -- (0,0.2);
 \node[below] at (0,-0.2) {$a = x_1$};

 \draw[thick] (4,-0.1) -- (4,0.1);
 \node[below] at (4,-0.2) {$x_2$};

 \draw[thick] (6,-0.1) -- (6,0.1);
 \node[below] at (6,-0.2) {$x_3$};

 \draw[thick] (7,-0.1) -- (7,0.1);
 \node[below] at (7,-0.2) {$x_4$};

 \node[below] at (8,-0.2) {$\dots$};

 \draw[thick] (9,-0.1) -- (9,0.1);
 \node[below] at (9,-0.2) {$x_n$};

 \draw[thick] (10,-0.2) -- (10,0.2);
 \node[below] at (10,-0.2) {$b$};
 \end{tikzpicture}
\end{center}
Note that, in this case, $\operatorname{C}_g([a,b],\mathbb F)$ is isomorphic to the space of convergent sequences. The following function clearly belongs to the closure of ${\rm P}_g$,
\[
\operatorname{exp}_g(-\Delta g(a)^{-1},a)(x):=\sum_{n=0}^\infty \frac{(-\Delta g(a))^{-n}}{n!}g_{a,n}(x)=\begin{dcases}
 1,&x=a,\\
 0,&x\in (a,b],
\end{dcases}
\]
see \cite[Proposition 5.7]{Cora2023}. Note that the function above is the indicator function of the set $\{a\}$. Furthermore, in the spirit of Theorem~\ref{equivalenciagpolinomiosdensidad}, the function
\[ g^*(x)=g_{a,1}(x)-\Delta g(a)(1-1_{\{a\}}(x))=\sum_{t\in D_g(a,t)}\Delta g(t)\]
is a derivator. We can approximate $g^*$ uniformly by $g$--polynomials so, we can approximate $g^*$-polynomials by $g$--polynomials since
\[
\int_{[a,x)} f\operatorname{d}\mu_{g^*} = \int_{[a,x)} f-f\cdot 1_{\{a\}}\operatorname{d}\mu_{g}.
\]
Hence, the function,
\[
\operatorname{exp}_{g^*}(-\Delta g(x_2)^{-1},a)(x)=\sum_{n=0}^\infty \frac{(-\Delta g(x_2))^{-n}}{n!}g^*_{a,n}(x)=\begin{dcases}
 1,&x\in [a,x_2],\\
 0,&x\in (x_2,b],
\end{dcases}
\]
belongs to the closure of ${\rm P}_g$. Subtracting the functions above we get
\[
1_{(a,x_2]}(x)=\begin{dcases}
 1,&x\in (a,x_2],\\
 0,&x\in \{a\}\cup(x_2,b].
\end{dcases}
\]
Proceeding by induction, eventually we get that every function of the form
\[
1_{(x_i,x_{i+1}]}(x)=\begin{dcases}
 1,&x\in (x_i,x_{i+1}],\\
 0,&x\notin (x_i,x_{i+1}],
\end{dcases}
\]
belongs to the closure of $\mathrm{P}_g$. These functions are clearly dense, so the set $\mathrm{P}_g$ must be dense.

\section{Final conclusions}
\label{finalcooncl}

In this work, we have extended the classical Weierstrass Approximation Theorem to the setting of Stieltjes calculus. Specifically, we have shown that when the derivator $g$ has only finitely many discontinuities, the space of $g$--polynomials is dense in the Banach space of uniformly $g$--continuous functions.

Our approach exploits the topological structure induced on the interval $[a,b]$ by a derivator with finitely many jumps. As established in Lemma~\ref{psi}, each discontinuity of $g$ disconnects the interval, effectively partitioning it into subintervals on which uniformly $g$--continuous functions—being regulated—are continuous. In this piecewise continuous setting, we showed how $g$--polynomials factor across the partition (Lemma~\ref{pgfactor}) and leveraged the classical Weierstrass theorem within each subinterval.

The argument crucially depends on the finiteness of the discontinuities. In the case of infinitely many jumps, the induced topology becomes substantially more complex, and the decomposition used in Lemma~\ref{psi} is no longer available. Nevertheless, this limitation does not preclude the validity of the theorem in the general case. In Section~\ref{secfinal}, we present a partial extension under topological assumptions on the discontinuity set, and in Section~\ref{seccionrutas}, we discuss the various proof strategies we explored, outlining a potential path for future work on the general case.

\section*{Funding}
The authors were partially supported by “ERDF A way of making Europe” of the “European Union”; and by Xunta de Galicia, Spain, project ED431C 2023/12 and grant PID2020-113275GB-I00 funded by MCIN/AEI/10.13039/501100011033, Spain. Victor Cora was partially supported by Xunta de Galicia under grant ED481A-2024-152.

\bibliographystyle{spmpsciper}
\bibliography{weierfinito}

\end{document}